\documentclass{article}[11pt]

\usepackage[ruled,vlined]{algorithm2e}

\usepackage{graphicx}
\usepackage{amsmath,amsfonts,amssymb,latexsym}
\usepackage{enumerate}
\usepackage{setspace}
\usepackage{color}
\usepackage[cp1250]{inputenc}
\usepackage{tikz}
\usepackage{soul}

\newtheorem{thm}{Theorem}[section]
\newtheorem{prop}[thm]{Proposition}

\newtheorem{cor}[thm]{Corollary}

\newtheorem{lemma}[thm]{Lemma}

\newtheorem{prob}{Problem}

\newcommand{\qed}{$\Box$}
\newcommand{\proof}{\noindent\textbf{Proof. }}

\newcommand{\cp}{\,\square\,}

\newcommand{\smallqed}{{\tiny ($\Box$)}}

\newenvironment{unnumbered}[1]{\trivlist \item [\hskip \labelsep {\bf
#1}]\ignorespaces\it}{\endtrivlist}

\def\vertex(#1){\put(#1){\circle*{1.8}}}
\def\lab(#1)#2{\put(#1){\makebox(0,0)[c]{#2}}}

\newcommand{\gr}{\gamma_{\rm gr}}

\newcommand{\grz}{\gamma_{\rm gr}^{\rm Z}}

\begin{document}

\title{Grundy domination and zero forcing in regular graphs}

\author{
Bo\v{s}tjan Bre\v{s}ar$^{a,b}$
\and
Simon Brezovnik$^{a,c}$
}

\date{\today}

\maketitle

\begin{center}
$^a$ Faculty of Natural Sciences and Mathematics, University of Maribor, Slovenia (bostjan.bresar@um.si, simon.brezovnik2@um.si)\\
$^b$ Institute of Mathematics, Physics and Mechanics, Ljubljana, Slovenia\\
$^c$ Faculty of Education, University of Maribor, Slovenia \\
\end{center}

\begin{abstract}
Given a finite graph $G$, the maximum length of a sequence $(v_1,\ldots,v_k)$ of vertices in $G$ such that each $v_i$ dominates a vertex that is not dominated by any vertex in $\{v_1,\ldots,v_{i-1}\}$ is called the Grundy domination number, $\gamma_{\rm gr}(G)$, of $G$. A small modification of the definition yields the Z-Grundy domination number, which is the dual invariant of the well-known zero forcing number. In this paper, we prove that $\gamma_{\rm gr}(G) \ge  \frac{n + \lceil \frac{k}{2} \rceil - 2}{k-1}$ holds for every connected $k$-regular graph of order $n$ different from $K_{k+1}$ and $\overline{2C_4}$. The bound in the case $k=3$ reduces to $\gr(G)\ge \frac{n}{2}$, and we characterize the connected cubic graphs with $\gamma_{\rm gr}(G)=\frac{n}{2}$. If $G$ is different from $K_4$ and $K_{3,3}$, then $\frac{n}{2}$ is also an upper bound for the zero forcing number of a connected cubic graph, and we characterize the connected cubic graphs attaining this bound.  
\end{abstract}



\section{Introduction}
\label{sec:intro}

Zero forcing is defined by the following process that starts by choosing a set $S$ of vertices of a graph $G$ and all vertices of $S$ are colored blue. Then the color-change operation is performed, in which a vertex, which is the only non-blue neighbor of a blue vertex, is colored blue. The color-change operation is performed as long as possible. If at the end of the process all vertices become blue, then the initial set $S$ is called a {\em zero forcing set} of $G$. The minimum cardinality of a zero forcing set in $G$ is the {\em zero forcing number}, $Z(G)$, of $G$. 
This concept arose in the study of minimum rank 
among symmetric matrices described by a graph~\cite{AIM}, and was rediscovered independently in mathematical physics and in graph search algorithms; see~\cite{BFF-2018} and the references therein. Zero forcing is closely  related to power domination, which was  introduced in~\cite{haynes-2002} as a model for monitoring electrical networks; cf.~\cite{BFF-2018}. It is also related to path-width and tree-width parameters~\cite{barioli, BBF2013,taklimi}, and recently its relation with the inverse eigenvalue problem was presented~\cite{fallat}. 

In the last decade, several Grundy domination invariants were introduced~\cite{bbgkkptv-2017, bgmrr-2014, bhr}, which were motivated by the domination games introduced in~\cite{brklra-2010, hkr}. An additional motivation for Grundy domination comes from the process of expanding a dominating set in a graph that is built on-line~\cite{bgk-2016}. The process can be presented by a sequence of vertices in a graph, and the basic version of Grundy domination from~\cite{bgmrr-2014} is defined as follows. A sequence $S=(v_1,\ldots,v_k)$ of vertices in a graph $G$ is a {\em closed neighborhood sequence} if for every $i\in \{2,\ldots,k\}$, 
\begin{equation}
\label{eq:defGrundy}
N_G[v_i] \setminus \bigcup_{j=1}^{i-1}N_G[v_j]\not=\emptyset,
\end{equation}
where $N_G[v_j]$ is the {\em closed neighborhood} of $v_j$ (note that $N_G[v_j]=N_G(v_j)\cup\{v_j\}$, where $N_G(v_j)$ is the set of neighbors of $v_j$). The corresponding set of vertices from the closed neighborhood sequence $S$ will be denoted by $\widehat{S}$. The maximum length $|\widehat{S}|$ of a closed neighborhood sequence $S$ in a graph $G$ is the {\em Grundy domination number}, $\gr(G)$, of $G$. Every maximal sequence $S$ enjoying the property~\eqref{eq:defGrundy} for all of its vertices is a {\em dominating sequence} in $G$. A vertex $x\in V(G)$ {\em dominates} a vertex $y$ if $y\in N_G[x]$, and we then also say that $y$ is {\em dominated} by $x$. If $D\subset V(G)$, then $y\in V(G)$ is {\em dominated} by $D$ if there exists $x\in D$ that dominates $y$. A set $D$ is a {\em dominating set} of a graph $G$ if every vertex in $G$ is dominated by $D$. 
Note that a closed neighborhood sequence $S$ is a dominating sequence in $G$ if and only if $\widehat{S}$ is a dominating set of $G$.
If $(v_1,\ldots,v_k)$ is a closed neighborhood sequence, then we say that $v_i$ \emph{footprints} the vertices from $N_G[v_i] \setminus \cup_{j=1}^{i-1}N_G[v_j]$, and that $v_i$ is the \emph{footprinter} of every vertex $u\in N_G[v_i] \setminus \cup_{j=1}^{i-1}N_G[v_j]$, for any $i\in [k]$ (where $[k]=\{1,\ldots,k\}$). 

Now, let $G$ be a graph with no isolated vertices. 
A closed neighborhood sequence $S$ in $G$ is a {\em Z-sequence} if, in addition, every vertex $v_i$ in $S$ footprints a vertex distinct from itself. Equivalently, $S$ is a Z-sequence in $G$ if 
$$N_G(v_i) \setminus \bigcup_{j=1}^{i-1}N_G[v_j]\not=\emptyset$$
holds for for every $i\in \{2,\ldots,k\}$. The maximum length of a Z-sequence in $G$ is the {\em Z-Grundy domination number}, $\grz(G)$, of $G$. Note that if $S$ is a Z-sequence, then saying that $x\in \widehat{S}$ {\em footprints} a vertex $y$ necessarily implies that vertices $x$ and $y$ are distinct.
The Z-Grundy domination number was introduced in~\cite{bbgkkptv-2017} as the dual of the zero forcing number.
Notably, $S$ is a Z-sequence if and only if the set of vertices outside $S$ forms a zero forcing set~\cite{bbgkkptv-2017}. In particular, 
\begin{equation}
\label{eq:zeroZGrundy}Z(G)=n(G)-\grz(G) 
\end{equation}
for every graph $G$ with no isolated vertices, where $n(G)$ is the order of $G$.  
In a subsequent paper, Lin presented a natural connection between four variants of Grundy domination and four variants of zero forcing~\cite{lin-2019}.  The connections show that all versions of Grundy domination can be applied in the study of different types of minimum rank parameters of symmetric matrices.  

One of the central problems concerning (domination) invariants is to find general bounds in terms of the order of a graph, possibly involving also the maximum degree or some other parameter. Interestingly, a general lower bound for the total version of the Grundy domination number of regular graphs was presented in~\cite{bhr}, but for the standard Grundy domination number such a bound has not yet been known. 
On the other hand, the Grundy domination number (and in some cases also its Z-variant) was studied in graph products~\cite{bbgkkptv-2016, nt-2018}, Sierpi\'{n}ski graphs~\cite{bgk-2016}, and in Kneser graphs~\cite{bkt-2019}. 

Several authors considered bounds on the zero forcing number in terms of the order, from which one can directly get dual bounds for the Z-Grundy domination number by using~\eqref{eq:zeroZGrundy}. Amos et al.~\cite{acdp-2015} proved the general upper bound, $Z(G)\le \frac{(\Delta-2)n+2}{\Delta-1}$, which holds for all connected graphs $G$ with maximum degree $\Delta$ and order $n$. 
Gentner et al.~\cite{gprs-2016} characterized the extremal graphs attaining the bound. Moreover, Gentner and Rautebach~\cite{gr-2018} improved the bound to $Z(G)\le \frac{n(\Delta - 2)}{\Delta - 1}$, whenever $G$ is a connected graph with $\Delta\ge 3$ that is not isomorphic to one of the five sporadic graphs presented by the authors. In particular, for graphs with $\Delta=3$, the bound reduces to $Z(G)\le \frac{n}{2}$. 
Gir\"{a}o et al.~\cite{gms-2018} presented an infinite family of graphs $G_n$ with maximum degree $3$ such that the zero forcing number of $G_n$ is at least $\frac{4}{9}n$. Davila and Henning studied the zero forcing number of connected claw-free cubic graphs $G$ and proved that $Z(G)<\frac{n}{2}$ as soon as $G$ has at least $10$ vertices~\cite{dh-2018}, and further improved this result in~\cite{dh-2020} to $Z(G)\le \frac{n}{3}+1$ by additionally excluding the 2-necklace graph (see the right graph $N_{YY}$ in Figure~\ref{fig:necklace}). 

In the next section, we will prove the following lower bound for the Grundy domination number of connected $k$-regular graphs of order $n$ different from $K_{k+1}$ and $\overline{2C_4}$: $$\gr(G) \ge  \frac{n + \lceil \frac{k}{2} \rceil - 2}{k-1}.$$
The result is similar to the bound from~\cite{bhr} for the total version of the Grundy domination number, and fills the gap in the study of the Grundy domination number. In Section~\ref{sec:ZGrundy} we prove a lower bound for the Z-Grundy domination number of regular graphs, which is similar to the bound from Amos et al.~\cite{acdp-2015}, yet it slightly improves it when $G$ has a triangle. Then, for connected cubic graphs $G$ different from $K_4$ and $K_{3,3}$ we prove that $\grz(G)\ge \frac{n}{2}$, which rediscovers (with a simpler proof) the bound of Gentner and Rautebach in~\cite{gr-2018} restricted to cubic graphs. Section~\ref{sec:cubic} contains our main result, which is a characterization of connected cubic graphs $G$ with $\grz(G)=\frac{n}{2}$. This result can be viewed as an extension of the results of Davila and Henning~\cite{dh-2018,dh-2020} from connected claw-free cubic to all connected cubic graphs with $Z(G)=\frac{n}{2}$. As a by-product, we also get the family of connected cubic graphs with $\gr(G)=\frac{n}{2}$. The extremal family in the former case contains $15$ sporadic graphs (see Figures~\ref{fig:X2X3}-\ref{fig:XY2} and \ref{fig:hamming}-\ref{fig:bip}), while in the later case reduces to $8$ graphs. 


\section{Grundy domination in regular graphs}
\label{sec:regular}

In this section, we establish a lower bound on the Grundy domination number of a regular graph. The Grundy domination number of cycles can be easily established, namely, $\gr(C_n)=n-2$. Hence, we may restrict to $k$-regular graphs with $k>2$. 

First, a few a more definitions. Vertices $u$ and $v$ in $G$ are {\em twins} if $N_G[u]=N_G[v]$ and are {\em open twins} if $N_G(u)=N_G(v)$. We write $kG$ for the disjoint union of $k$ copies of a graph $G$. The complement of a graph $G$ is denoted by $\overline{G}$.

\begin{thm}
If $k \ge 3$ and $G$ is a connected $k$-regular graph of order~$n$ different from $K_{k+1}$ and $\overline{2C_4}$, then
$$\gr(G) \ge  \frac{n + \lceil \frac{k}{2} \rceil - 2}{k-1}. $$
\label{thm:regular}
\end{thm}
\proof Let $k\ge 3$ and let $G$ be a connected $k$-regular graph of order~$n$ different from $K_{k+1}$ and $\overline{2C_4}$. 
Suppose that $G$ is bipartite, and let $V(G)=A\cup B$ be the bipartition of the vertex set of $G$ into independent sets $A$ and $B$. Since $\gr(G)\ge \alpha(G)\ge \frac{n}{2}$, we infer 
$$\gr(G)\ge \frac{n}{2}=\frac{n+\lceil\frac{3}{2}\rceil-2}{3-1}\ge\frac{n + \lceil \frac{k}{2} \rceil - 2}{k-1},$$
\noindent which is true for any $k\ge 3$. Thus, we may assume that $G$ is non-bipartite. 

It is clear that $\gr(G)>1$, because $G$ is not a complete graph. In addition, $\gamma(G)>1$, for otherwise $G$ is not regular (since it is not complete). Suppose that $\gr(G)=2$.
By~\cite[Theorem 3.6]{bgmrr-2014}, $\gr(G)=2=\gamma(G)$ if and only if $G$ is the join of one or more graphs $\overline{K_{r,s}}$. Since $G$ is regular, all of these graphs are $\overline{K_{r,r}}$ for some fixed integer $r$. So let $G$ be the join of $\ell$ graphs $\overline{K_{r,r}}$.  Since $G$ is connected, $\ell\ge 2$.  Note that $k=(\ell-1)2r+r-1$, and $n=2r\ell$, which implies
$$\frac{n + \lceil \frac{k}{2} \rceil - 2}{k-1}= \frac{3r\ell-r+\lceil\frac{r-1}{2}\rceil-2}{2r\ell-r-2}.$$
\noindent Now, since $\ell\ge 2$, we have $$\frac{3r\ell-r+\lceil\frac{r-1}{2}\rceil-2}{2r\ell-r-2}\le 2=\gr(G)$$
\noindent if and only if $(r,\ell)\notin\{(1,2),(2,2)\}$. The first case, $r=1,\ell=2$, gives $G=C_4$, which is $2$-regular and not relevant for this proof. The second case, $r=2,\ell=2$, gives $G=\overline{2C_4}$, which is also excluded in the assumption of the statement of the theorem. In the rest of the proof we may thus assume that $\gr(G)>2$. We distinguish two cases. 

{\bf Case 1.} $G$ has a triangle.

It is clear that there exist two vertices that lie in a triangle in $G$ that are not twins. (Indeed, if every two vertices that lie in a triangle in $G$ were twins, then since $G$ is connected, this would imply that $G$ is complete, a contradiction.) Let $v_1$ and $v_2$ be chosen among all adjacent vertices in $G$ that are not twins to have the maximum number of common neighbors (as noted above, they have at least one common neighbor). Clearly, there are at most $k-2$ common neighbors of $v_1$ and $v_2$, since they are not twins. We build the sequence $S$ starting with $(v_1,v_2)$. Note that $v_2$ footprints at most $k-2$ vertices, since $v_1$ and $v_2$ are already dominated when $v_2$ is added to the sequence. After the $i$th vertex is added to $S$, the sequence is $(v_1,\ldots,v_i)$, where $i\ge 2$. Suppose that this is not yet the entire sequence $S$, that is, $\{v_1,\ldots,v_i\}$ is not a dominating set of $G$. We choose $v_{i+1}$ as a vertex, which footprints at least one, but a minimum number of vertices. We claim that such a vertex $v_{i+1}$ exists. 

Let $x$ be any vertex not dominated by $\{v_1,\ldots,v_{i}\}$. Since $G$ is connected, there exists a path from a vertex $v_j$, where $j\in [i]$, to $x$, and consider a shortest possible such path $P$. Then the neighbor $y$ of $v_j$ on $P$ has the desired property, since it footprints the neighbor on $P$ different from $v_j$ (this is because $P$ is chosen as the shortest possible path between vertices of $\{v_1,\ldots,v_{i}\}$ and $x$). Hence, such a vertex $y$ exists, and footprints at most $k-1$ vertices. Thus $v_{i+1}$ is well defined and footprints at most $k-1$ vertices. 

By using the above construction, let $S=(v_1,\ldots,v_t)$ be the resulting dominating sequence of $G$. 

\begin{unnumbered}{Claim~A.}
One of the vertices $v_2$ or $v_t$ footprints at most $\lfloor\frac{k}{2}\rfloor$ vertices.  
\end{unnumbered}
\textbf{Proof of Claim~A.} Suppose to the contrary that each of the vertices $v_2$ and $v_t$ footprints at least $\lfloor\frac{k}{2}\rfloor+1$ vertices. In particular, this implies that $v_1$ and $v_2$ have at most $\lfloor\frac{k}{2}\rfloor-1$ common neighbors (recall that, by the choice of $v_1$ and $v_2$, this is also the maximum number of common neighbors that two non-twin neighbors in $G$ may have).  

Let $F$ be the set of vertices footprinted by $v_t$, and let $U=V(G)\setminus F$. First, suppose that $v_t\notin F$; that is, $v_t$ has been dominated by vertices in $\{v_1,\ldots,v_{t-1}\}$. Let $v_j$, $j\in [t-1]$, be the vertex that footprints $v_t$. We claim that $F$ induces a complete graph. Suppose that there exists $x\in F$ that is not adjacent to all other vertices in $F$. Then $x$ footprints less vertices than $v_t$, which is a contradition, since $x$ would be a better choice than $v_t$ for adding to the sequence $S$. Now, vertices $x$ and $v_t$ are neighbors, which are not twins (since $xv_j\notin E(G)$ and $v_tv_j\in E(G)$), and have $|F|-1\ge \lfloor\frac{k}{2}\rfloor$ common neighbors. This implies that $x$ and $v_t$ are non-twin neighbors with more common neighbors as $v_1$ and $v_2$, which is a contradiction to the choice of $v_1$ and $v_2$. 

Now, suppose that $v_t\in F$. Since $F\neq V(G)$ and $G$ is connected, there is a vertex $y\in U$, which is adjacent to a vertex in $F$. Note that $y$ is adjacent to all vertices of $F$, for otherwise $y$ would be a better choice than $v_t$ to be added to $S$. Hence, instead of $v_t$ we put $y$ as the last vertex of $S$, and we are in the situation of the previous paragraph, where $v_t\notin F$. In either case, the assumption that each of the vertices $v_2$ and $v_t$ footprints at least $\lfloor\frac{k}{2}\rfloor+1$ vertices leads us to a contradiction. 
~\smallqed

The next claim is about the case $\gr(G)\ge t>3$, however, along the way we will verify also the situation when $t=3$. 

\begin{unnumbered}{Claim~B.}
If $t>3$, then the sum of the numbers of vertices footprinted by $v_2,v_3$ and $v_t$ is at most $2k-3+\lfloor\frac{k}{2}\rfloor$.  
\end{unnumbered}
\textbf{Proof of Claim~B.}
As noted in the beginning of the construction of the sequence $S$, vertex $v_2$ footprints at most $k-2$ vertices. By Claim~A, $v_2$ or $v_t$ footprints at most $\lfloor\frac{k}{2}\rfloor$ vertices. If this is true for $v_t$, then the claim is proven, since $v_3$ footprints at most $k-1$ vertices, which was proved to hold for all vertices of $S$ except $v_1$. It remains to consider the possibility when $v_2$ footprints exactly $\lfloor\frac{k}{2}\rfloor$ vertices and $v_t$ footprints exactly $k-1$ vertices. 

Let $A=N(v_1)\setminus N[v_2]$, $C=N(v_2)\setminus N[v_1]$, and let $B$ be the set of common neighbors of $v_1$ and $v_2$. Since $t>2$ and $G$ is connected, there exists $y\in V(G)\setminus (N[v_1]\cup N[v_2])$, which is adjacent to a vertex $x$ in $A\cup B\cup C$. If $x\in B$, then we let $v_3=x$, and note that $x$ does not footprint $v_1, v_2$, and itself, therefore it footprints at most $k-2$ vertices. Thus, $v_2,v_3$ and $v_t$ footprint at most $(\lfloor\frac{k}{2}\rfloor)+(k-2)+(k-1)$ vertices, as claimed. (Note that if $t=3$, then $n\le (k+1)+(\lfloor\frac{k}{2}\rfloor)+(k-2)$, which implies $\frac{n + \lceil \frac{k}{2} \rceil - 2}{k-1}\le \frac{k+1+\lfloor\frac{k}{2}\rfloor+k-2+\lceil\frac{k}{2}\rceil-2}{k-1}=3\le\gr(G)$, in which case the statement of the theorem is correct.)

We may thus assume that there is no such vertex $x\in B$, and so all vertices in $B$ are adjacent only to vertices in $N[v_1]\cup N[v_2]$. 

Suppose that $x\in A\cup C$, and without loss of generality let $x\in A$. Let $|A|=s=|C|$. Note that $|B|=k-s-1$. If $B$ induces a complete graph, then from a vertex $b\in B$, there are $k-2-(k-s-2)=s$ edges connecting $b$ to vertices in $A\cup C$. In any case, even if $B$ is not a clique, there are at least $s$ edges connecting each vertex $b\in B$ to vertices in $A\cup C$. Now, if there exists $b\in B$ such that all $s$ edges from $b$ to $A\cup C$ lead to vertices of $A$, then $b$ is adjacent also to $x$ (since $|A|=s$). Then, $v_3=x$, in which case $v_3$ footprints $y$, but it footprints at most $k-2$ vertices, which again proves the claim. (Again, in the case $t=3$, we derive the same inequality as earlier, confirming the correctness of the theorem.)

Finally, assume that there is an edge $e$ from $B$ to $C$, and let $z\in C$ be its endvertex. If $zx\in E(G)$, then we conclude the proof as in the previous paragraph, by taking $v_3=x$. Hence, let us assume that $z$ is adjacent to at most $s-1$ vertices in $A$ (that is, $zx\notin E(G)$). We may also assume that $z$ does not have a neighbor in $V(G)\setminus (N[v_1]\cup N[v_2])$, because otherwise, by letting $v_3=z$ we derive the same conclusion as before, since then $z$ would footprint at most $k-2$ vertices. Hence, $z$ has $k-1-(s-1)=k-s$ neighbors in $B\cup C$. Note that $v_2$ and $z$ are neighbors, which are not twins (since $v_1v_2\in E(G)$ and $v_1z\notin E(G)$), and have $k-s$ common neighbors. Since $v_1$ and $v_2$ have only $k-s-1$ common neighbors, this is a contradiction with the choice of $v_1$ and $v_2$. This implies the truth of the claim. 
~\smallqed
\medskip

Note that in the proof of the above claim, we have also proved the statement of the theorem for the case $\gr(G)=t=3$ (when $G$ has a triangle), hence we may assume $t>3$. 
Since $v_1$ footprints $k+1$ vertices, and, by Claim~B, $v_2,v_3$ and $v_t$ together footprint at most $2k-3+\left\lfloor\frac{k}{2}\right\rfloor$ vertices, and all other $t-4$ vertices each footprints at most $k-1$ vertices, we infer:
\[
n\le (k+1)+(2k-3+\lfloor\frac{k}{2}\rfloor)+(t-4)(k-1),
\]

\noindent which gives

\[
\gr(G)\ge t \ge \frac{n+k-\lfloor\frac{k}{2}\rfloor-2}{k-1}=\frac{n + \lceil \frac{k}{2} \rceil - 2}{k-1}. 
\]

{\bf Case 2.} $G$ is triangle-free.

The basic part of the construction of a dominating sequence $S=(v_1,\ldots,v_t)$ is similar as in Case 1. Note that, since $G$ has no triangles, no two adjacent vertices are twins. If $i<t$, and we have constructed $(v_1,\ldots,v_i)$, then we choose $v_{i+1}$ as a vertex, which footprints at least one, but a minimum number of vertices. In the same way as in Case 1 we can prove that such a vertex $v_{i+1}$ exists and footprints at most $k-1$ vertices. 

The construction of $S$ is based on the rule that each vertex added to the sequence footprints at least one, but the least possible number of vertices. 
Clearly, $v_1$ footprints $k+1$ vertices, and as noted above, each further vertex in $S$ footprints at most $k-1$ vertices. Next, we will prove that the last vertex of the sequence footprints just one vertex. 

\begin{unnumbered}{Claim~C.}
Vertex $v_t$ footprints only one vertex.  
\end{unnumbered}
\textbf{Proof of Claim~C.}
Suppose to the contrary that $v_t$ footprints two distinct vertices $x$ and $y$. Let us first assume that $v_t$ is dominated before it is added to the sequence. Clearly, $x$ and $y$ are not adjacent, since $G$ has no triangles. But then $x$ would be a better choice than $v_t$ for adding it to $(v_1,\ldots,v_{t-1})$, since it footprints less vertices than $v_t$, which is a contradiction. Second case is that $v_t$ is not dominated by vertices of $\{v_1,\ldots,v_{t-1}\}$. Hence, $v_t$ footprints itself, and we may write $x=v_t$, and so $v_t$ footprints also $y$. Since $G$ is connected, there is a vertex $u$, adjacent to $v_t$, which is dominated by $\{v_1,\ldots,v_{t-1})$. If $uy\notin E(G)$, then again $u$ is a better choice than $v_t$ to be added to $(v_1,\ldots,v_{t-1})$, which contradicts the construction. Thus, $u$ is adjacent to both $x$ and $y$, and we can put $u$ at the end of the sequence $S$ instead of $v_t$. But then we are in the previous case, where $v_t$ does not footprint itself. 
~\smallqed

Suppose that $k\ge 4$. Then $1\le \lfloor\frac{k}{2}\rfloor-1$. Summing the upper bounds on the number of vertices, which are footprinted at each step, we infer 
$$n\le k+1+(t-2)(k-1)+1\le k+1+(t-2)(k-1)+\lfloor\frac{k}{2}\rfloor-1,$$
\noindent
which implies
$$\gr(G)\ge t \ge  \frac{n + \lceil \frac{k}{2} \rceil - 2}{k-1}. $$

Finally, we are left with the case when $G$ is a $3$-regular (triangle-free) graph. Since $G$ is non-bipartite, there exist odd cycles, and let $C$ be a shortest odd cycle in $G$. This implies, in particular, that $C$ is induced. Since $G$ is connected and $3$-regular, there are vertices outside $C$ that are adjacent to vertices in $C$. 

Suppose that there exists a vertex $x\in V(G)$, which is not on $C$, and is adjacent to two vertices of $C$, say $u$ and $v$. Among all such possible vertices let $x$ be chosen in such a way that two neighbors $u$ and $v$ on $C$ are as close as possible with respect to the distance on $C$. Since $G$ has no triangles, $u$ and $v$ are not adjacent. Let $u=v_1,v_2,\ldots,v_r=v$ be a shortest path on $C$ between $u$ and $v$. We start the sequence $S$ in the same way, that is, with $(v_1,\ldots, v_r)$. Clearly, each of the vertices $v_i$, where $1\le i\le r$ footprints $v_{i+1}$, because $C$ is an induced cycle. Moreover, $v_r=v$ footprints only one vertex,  namely, $v_{r+1}$, since $x$ has already been footprinted by $v_1=u$.

The second case is that every vertex $x\in V(G)$, which is not on $C$ has at most one neigbbor on $C$. Therefore, each vertex of $C$ has a unique neighbor that is not on $C$. Therefore, letting $C:v_1,\ldots,v_s,v_1$, where vertices are written in the natural order, we can start the sequence $S$ with $(v_1,\ldots,v_s)$. Clearly, $v_s$ footprints only one vertex, the neighbor of $v_s$ outside $C$. Now, if $v_s$ is the last vertex of the sequence $S$, then $n=2s$, and we derive that $$\gr(G)\ge s=\frac{n}{2}=\frac{n + \lceil \frac{3}{2} \rceil - 2}{3-1},$$ 
as desired. Otherwise, the sequence $S$ has more vertices, that is, $t>s$.

In either of the above two cases (last two paragraphs), we found a vertex in $S$, which is not $v_t$, that footprints only one vertex. Noting that $v_1$ footprints $k+1=4$ vertices, $v_t$ footprints $1$ vertex, and all other vertices footprint at most $k-1=2$ vertices, we infer 
$$n\le 4+1+1+2(t-3),$$
which implies 
$$\gr(G)\ge t\ge \frac{n}{2}=\frac{n + \lceil \frac{3}{2} \rceil - 2}{3-1}.$$ 
The proof is complete. ~\qed

\bigskip



In the special case of Theorem~\ref{thm:regular} when $k = 3$ have the following result.

\begin{cor}
If $G \ne K_{4}$ is a cubic graph of order~$n$, then $\gr(G) \ge \frac{1}{2}n$,
and the bound is sharp. \label{cor:cubic}
\end{cor}

The Petersen graph is a well-known example that attains the bound in Corollary~\ref{cor:cubic}. As we will see in Section~\ref{sec:cubic}, where we will characterize the cubic graphs $G$ with $\gr(G)=\frac{n}{2}$, there are $8$ such graphs.

The question is, whether the bound of Theorem~\ref{thm:regular} is sharp also for $k\ge 4$, and to characterize all extremal graphs. 

\begin{prob}
If $k \ge 4$, determine all  $k$-regular graphs $G$ of order~$n$ such that
$$\gr(G)=\frac{n + \lceil \frac{k}{2}\rceil - 2}{k-1}. $$
\end{prob}


\section{Z-Grundy domination in regular graphs}
\label{sec:ZGrundy}

In a similar way as in Theorem~\ref{thm:regular}, we can prove the following bound for the Z-Grundy domination number of regular graphs.

\begin{thm}
If $k \ge 3$ and $G$ is a connected $k$-regular graph of order~$n$ different from $K_{k+1}$, then
$$\grz(G) \ge  \left\{
        \begin{array}{cl}
                 \frac{n-1}{k-1} ;& G \textrm{ has a triangle,}  \\ 
                           
                \frac{n-2}{k-1}; & G \textrm{ triangle-free}
         \end{array}
       \right..$$
\label{thm:ZGrundy}
\end{thm}
\proof Let $k\ge 3$ and let $G$ be a connected $k$-regular graph of order~$n$ different from $K_{k+1}$.  We distinguish two cases, depending of whether $G$ is triangle-free.

First, let $G$ have a triangle. 
As noted in the proof of Theorem \ref{thm:regular}, there exist two vertices that lie in a triangle in $G$ that are not twins. Let $v_1$ and $v_2$ be chosen among all adjacent vertices in $G$ that are not twins to have the maximum number of common neighbors. As noted above, they have at least one common neighbor, hence $|N_G[v_2]\setminus N_G[v_1]|\le k-2$.
If $\{v_1,v_2\}$ is a dominating set of $G$, then $n\le k+1+k-2$, which implies $$\grz(G)\ge 2\ge \frac{n-1}{k-1}.$$

We may thus assume that $\grz(G)>2$, and we build a Z-sequence $S$ starting with $(v_1,v_2)$. As noted above, $v_2$ footprints at most $k-2$ vertices. After the $i$th vertex is added to $S$, the sequence is $(v_1,\ldots,v_i)$, where $i\ge 2$. Suppose that this is not yet the entire sequence $S$, that is, $\{v_1,\ldots,v_i\}$ is not a dominating set of $G$. We choose $v_{i+1}$ as a vertex, which footprints at least one, but a minimum number of vertices. In the same way as in the proof of Theorem~\ref{thm:regular}, we can prove that such a vertex $v_{i+1}$ exists. In addition, since it can be chosen among vertices, that have been dominated by the set $\{v_1,\ldots,v_i\}$, we infer that $v_{i+1}$ footprints at most $k-1$ vertices. Let $t$ be the length of $S$. Then, $$n\le (k+1)+(k-2)+(t-2)(k-1),$$ which implies $$\grz(G)\ge t\ge \frac{n-1}{k-1}.$$

The second case, when $G$ is triangle-free can be proved almost in the same way. Note, however, that we cannot assume that $v_2$ footprints at most $k-2$ vertices. Yet, we find that it footprints at most $k-1$ vertices, because  $v_1$ and $v_2$ can be chosen as neighbors that are not twins (since $G$ is triangle-free and $k\ge 3$, every pair of adjacent vertices are good). Now, counting the number of vertices footprinted in each step, we infer 
$$n\le (k+1)+(t-1)(k-1),$$ which implies $$\grz(G)\ge t\ge \frac{n-2}{k-1}.$$
\qed 

\bigskip

The bound in Theorem~\ref{thm:ZGrundy} is sharp in the case of triangle-free graphs. Note that for every $k\ge 3$, we have $\grz(K_{k,k})=2=\frac{2k-2}{k-1}=\frac{n-2}{k-1}$.
When $G$ has a triangle, we note that $\grz(G)$ is an integer, and so we infer the slightly improved bound $\grz(G)\ge \lceil\frac{n-2}{k-1}\rceil$. In this case, $\grz(K_3\cp K_2)=3=\lceil\frac{6-1}{3-1}\rceil= \lceil\frac{n-1}{k-1}\rceil$.  

Applying~\eqref{eq:zeroZGrundy}, the following consequence for the zero forcing number of regular graphs is immediate.
\begin{cor}
If $k \ge 3$ and $G$ is a connected $k$-regular graph of order~$n$ different from $K_{k+1}$, then
$$Z(G) \le \left\{
        \begin{array}{cl}
                \frac{n(k-2)+1}{k-1} ;& G \textrm{ has a triangle,}  \\ 
                           
                \frac{n(k-2)+2}{k-1}; & G \textrm{ triangle-free}
         \end{array}
       \right..$$
\label{cor:Zforcing}
\end{cor}

The bound in the triangle-free case coincides with the bound of Amos et al.~\cite{acdp-2015} for connected graphs with maximum degree $\Delta$: 
$$Z(G)\le \frac{(\Delta-2)n+2}{\Delta-1}.$$
It was also proved in~\cite{gprs-2016} that the bound is attained if and only if $G$ is either $K_n$, $C_n$ or $K_{\Delta,\Delta}$. 

The lower bound on the Z-Grundy domination number of regular graphs obtained in Theorem~\ref{thm:ZGrundy} restricted to cubic graphs $G$ states $\grz(G)\ge \frac{n-2}{2}$. We next improve the bound for cubic graphs as soon as they are different from $K_4$ and $K_{3,3}$. For the proof of this result, we need the following lemma. 

\begin{lemma}
\label{lem:Zsequence}
Let $G$ be a connected cubic graph of order $n$. If there exists a  (not necessarily dominating) Z-sequence $S=(v_1,\ldots,v_r)$ in $G$ such that each vertex $v_i$, where $i\in \{2,\ldots, r\}$, footprints at most two vertices,  and there are at least two vertices in $S$ that footprint only one vertex, then $\grz(G)\ge\frac{n}{2}$. 
\end{lemma}
\proof Note that $v_1$ footprints $4$ vertices, and by the assumption of the lemma, vertices $v_1,\ldots,v_r$ together footprint at most $4+2+2(r-3)$ vertices. If $S$ is already a dominating sequence, then the proof can be continued in the last paragraph. Otherwise, we will prove that $S$ can be extended to a dominating sequence in such a way that for all the remaining vertices each footprints at most two vertices. 

Let $D$ be the set of vertices that are dominated by $S$, which is not a dominating sequence. Since $G$ is connected, it contains a vertex $x\notin D$, which is adjacent to a vertex $y\in D$. Clearly, $y\notin\widehat{S}$, and let $v_i$ be a vertex in $S$ that dominates $y$. If $y$ is added to $S$, that is, $v_{r+1}=y$, then $y$ footprints at most two vertices, because $\{v_i,y\}\subset N_G[y]$ and $|N_G[y]|=4$. Since $y$ footprints $x$, we may let $v_{r+1}=y$. Repeating this argument, we can easily see that $S$ can be extended to a dominating sequence of length $t$ in such a way that each of $v_{r+1},\ldots, v_t$ footprint at most two vertices.

In the same way as in the beginning of this proof we derive that vertices $v_1,\ldots,v_t$ together footprint at most $4+2+2(t-3)$ vertices, which implies $n\le 4+2+2(t-3)$, and we infer $$\grz(G)\ge t\ge \frac{n}{2}.$$
\qed

We are ready for the proof of the lower bound for the Z-Grundy domination number of cubic graphs (see also~\cite{gr-2018}, where a more general result is proved for the zero forcing number, but with a more difficult proof).

\begin{thm}
If $G$ is a connected cubic graph of order~$n$ different from $K_{4}$ and $K_{3,3}$, then $\grz(G) \ge   \frac{n}{2}$ and the bound is sharp. 
\label{thm:Z-cubic}
\end{thm}
\proof  Let $G$ be a connected cubic graph of order $n$ such that $G\ne K_4$ and $G\ne K_{3,3}$. In view of Lemma~\ref{lem:Zsequence} it suffices to find a Z-sequence $S$ (not necessarily dominating) in $G$ such that each vertex of $S$ except the first vertex footprints at most two vertices, and there are two vertices that footprint only one vertex. We distinguish two cases with respect to $G$ having a triangle or not. 

Suppose that $G$ has a triangle, but $G$ does not have a diamond. Let a triangle in $G$ have vertices $v_1,v_2,v_3$. Then each of $v_i,i\in [3]$, has a neighbor $u_i$, which is not a neighbor of $v_j$ for $j\ne i$. Hence, the sequence $(v_1,v_2,v_3)$ is a Z-sequence, and vertices $v_2$ and $v_3$ each footprints only one vertex, namely, $u_2$ and $u_3$, respectively. On the other hand, it is possible that $G$ has diamonds, and let $v_1,v_2,v_3,a$ be a diamond in $G$ such that $av_1\notin E(G)$ (we may assume this, since $G$ is not a complete graph). If $N_G(v_1)\ne N_G(a)$, then the sequence $(v_1,v_2,a)$ is a Z-sequence, and $v_2$ and $a$ each footprint only one vertex. Otherwise, let $N_G(v_1)=N_G(a)=\{v_2,v_3,x\}$. In this case, the sequence $(v_1,v_2,x)$ is a Z-sequence, and note that each of the vertices $v_2$ and $x$ footprint only one vertex.

Let $C:v_1,\ldots,v_r,v_1$ be a shortest cycle in $G$, and let $k\ge 4$. Suppose that for each $i\in [r]$ there exists $u_i\in N(v_i)$ such that $u_i$ is not a neighbor of $v_j$ for any $j\ne i$. Then $(v_1,\ldots,v_r)$ is a Z-sequence, since $v_i$ footprints $u_i$ for all $i\in [r]$. In addition, $v_{r-1}$ and $v_r$ each footprint only one vertex, and applying Lemma~\ref{lem:Zsequence}, we are done. Note that if $r\ge 5$, then the vertices $v_i,i\in [r]$ indeed have the property that each of their neighbors $u_i,i\in[r]$ outside $C$ is unique. Hence, we are left with the case when $C$ is a $4$-cycle. Since $G$ has no triangles it is not possible that $v_i$ and $v_{i+1}$ have a common neighbor (where $i$ is taken modulo $4$). Hence, it is only possible that $v_1$ and $v_3$ have a common neighbor and/or $v_2$ and $v_4$ have a common neighbor. Taking symmetry into account, there are two possibilities. 

First, if $v_1$ and $v_3$ have a common neighbor $b\notin V(C)$, but $N_G(v_2)=\{v_1,v_3,u_2\}$, $N_G(v_4)=\{v_1,v_3,u_4\}$, and $u_2\ne u_4$. Note that $b$ cannot be adjacent to both $u_2$ and $u_4$, and assume without loss of generality that $bu_2\notin E(G)$. If also $bu_4\notin E(G)$, then $(v_1,b,v_2,v_4)$ is a Z-sequence, in which $v_2$ and $v_4$ each footprints only one vertex. On the other hand, let $bu_4\in E(G)$. Now, if $u_2u_4\in E(G)$, then in the sequence $(v_1,b,u_4,u_2)$, each of the vertices $u_4$ and $u_2$ footprints only one vertex. Otherwise, if $u_2u_4\notin E(G)$, then in the sequence $(v_1,b,v_2,u_4)$, each of the vertices $v_2$ and $u_4$ footprints only one vertex. 

Second, let $v_1$ and $v_3$ have a common neighbor $b\notin V(C)$, and let  $v_2$ and $v_4$ have a common neighbor $c\notin V(C)$. Since $G\ne K_{3,3}$, $b$ and $c$ are not adjacent. Now, if $b$ and $c$ have no common neighbor, then in the sequence $(v_1,v_2,b,c)$, each of the vertices $b$ and $c$ footprints only one vertex. If, on the other hand, $b$ and $c$ have a common neighbor $x$, then in the sequence $(v_1,v_2,b,x)$, each of the vertices $b$ and $x$ footprints only one vertex. The proof is complete.
\qed

\bigskip

The sharpness of the bound in the above theorem is discussed in the next section, where we also characterize the graphs attaining the bound. 


\section{Cubic graphs attaining the $\frac{n}{2}$-bound}
\label{sec:cubic}

We start with the following auxiliary result, which is a slight modification of Lemma~\ref{lem:Zsequence}. 

\begin{lemma}
\label{lem:ZsequenceVec}
Let $G$ be a connected cubic graph of order $n$. If there exists a  (not necessarily dominating) Z-sequence $S=(v_1,\ldots,v_r)$ in $G$ such that each vertex $v_i$, where $i\in \{2,\ldots, r\}$, footprints at most two vertices,  and there are at least three vertices in $S$ that footprint only one vertex, then $\grz(G)>\frac{n}{2}$. 
\end{lemma}
\proof Note that $v_1$ footprints $4$ vertices, and by the assumption of the lemma, vertices $v_1,\ldots,v_r$ together footprint at most $4+3+2(r-3)$ vertices. If $S$ is already a dominating sequence, then the bound $\grz(G)>\frac{n}{2}$ holds, as argued in the last paragraph of this proof. Otherwise, we claim that $S$ can be extended to a dominating sequence in such a way that for all the remaining vertices each footprints at most two vertices. 

Let $D$ be the set of vertices that are dominated by $S$, where $S$ is a Z-sequence, but not a dominating sequence. Since $G$ is connected, it contains a vertex $x\notin D$, which is adjacent to a vertex $y\in D$. Clearly, $y\notin\widehat{S}$, and let $v_i$ be a vertex in $S$ that dominates $y$. If $y$ is added to $S$, that is, $v_{r+1}=y$, then $y$ footprints at most two vertices, because $\{v_i,y\}\subset N_G[y]$ and $|N_G[y]|=4$. Since $y$ footprints $x$, we may choose $v_{r+1}=y$. Repeating this argument, we can easily see that  $S$ can be extended to a dominating sequence of length $t$ in such a way that each of $v_{r+1},\ldots, v_t$ footprint at most two vertices.

In the same way as in the beginning of this proof we derive that vertices $v_1,\ldots,v_t$ together footprint at most $4+3+2(t-3)$ vertices, which implies $n\le 4+3+2(t-3)$, and we infer $$\gr(G)\ge t\ge \frac{n+1}{2}>\frac{n}{2}.$$
The proof is complete.  
\qed

\bigskip

The above lemma will be continuously used in the subsequent proofs in this section. In many cases we will determine a Z-sequence in which there are three vertices each of which footprints only one vertex while all other vertices, expect the first vertex of the sequence, footprint at most two vertices. By Lemma~\ref{lem:ZsequenceVec} this will imply that $\grz(G)>\frac{n}{2}$.

\begin{figure}[ht!]
\begin{center}
\begin{tikzpicture}[scale=0.52,style=thick]
\def\vr{4pt}

\path (6,4) coordinate (e'2);
\path (6,2) coordinate (a'2);
\path (8,1) coordinate (b'2);
\path (4,1) coordinate (c'2);
\path (6,0) coordinate (d'2);
\path (-6,2) coordinate (a');
\path (-6,0) coordinate (a'');
\path (-8,1) coordinate (b');
\path (-2,1) coordinate (c');
\path (-4,2) coordinate (d');
\path (-4,0) coordinate (d'');
\path (-5,4.5) coordinate (e);

\draw (6,4.6) node [text=black]{$\alpha_Y$}; 
\draw (-5,5.1) node [text=black]{$\alpha_X$};

\draw (6,-1.2) node [text=black]{\Large $Y$};
\draw (-5,-1.2) node [text=black]{\Large $X$};

\draw (a'2) -- (b'2) -- (d'2)--(c'2)--(a'2); 
\draw (a'2) -- (d'2);
\draw (b'2)--(e'2)--(c'2); 


\draw (d'')--(a') -- (b'); 
\draw (d')--(c');
\draw (a') -- (d')--(a''); 
\draw (b')--(a'')--(d'')--(c');
\draw (b')--(e)--(c');

\draw (a'2) [fill=white] circle (\vr);
\draw (b'2) [fill=white] circle (\vr);
\draw (c'2) [fill=white] circle (\vr);
\draw (d'2) [fill=white] circle (\vr);
\draw (e'2) [fill=black,] circle (\vr);


\draw (a') [fill=white] circle (\vr);
\draw (a'') [fill=white] circle (\vr);
\draw (b') [fill=white] circle (\vr);
\draw (c') [fill=white] circle (\vr);
\draw (d') [fill=white] circle (\vr);
\draw (d'') [fill=white] circle (\vr);
\draw (e) [fill=black] circle (\vr);

\end{tikzpicture}
\end{center}
\caption{Graphs $X$ and $Y$ with designated vertices marked.}
\label{fig:MandX}
\end{figure}
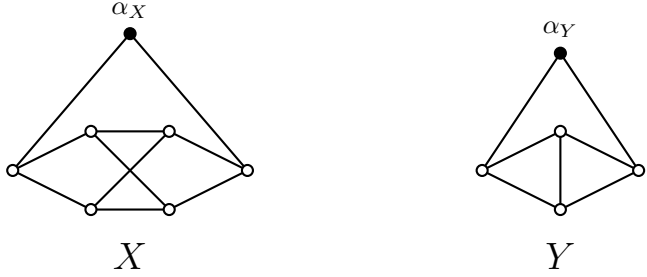

First, we will determine the graphs $G$ that attain the value $\grz(G)=\frac{n}{2}$ within a special family $\cal M$ of cubic graphs, which will be used in the proof of the main theorem. We need a couple of more definitions to introduce this family.

The graph $X$ is obtained from the complete bipartite graph $K_{3,3}$ by subdividing one edge. The vertex of degree $2$ in $X$ is denoted by $\alpha_X$. See Figure~\ref{fig:MandX}, where $X$ is depicted on the left. 
The graph $Y$ is obtained from the complete bipartite graph $K_{2,3}$ by adding an edge between two vertices of degree $2$; note that $Y$ has exactly one vertex of degree $2$, which we denote by $\alpha_Y$, while other four vertices of $Y$ have degree $3$.  One can find the graph $Y$ on the right side of Figure~\ref{fig:MandX}. 

The family $\cal M$ of (cubic) graphs is obtained from bi-regular trees $T$ with $\deg(u)\in \{1,3\}$ for every vertex $u\in V(T)$ in the following way. To every leaf of $T$ associate either a graph $X$ or a graph $Y$, and then attach a copy of the associated graph to that leaf. More precisely, $G\in \cal M$ if $G$ can be obtained from a tree $T$ whose non-leaf vertices have degree $3$ by identifying each leaf $\ell$ of $T$ with the vertex $\alpha_X$ or $\alpha_Y$ of a copy of $X$ or $Y$, respectively, that is associated with $\ell$.

\begin{figure}[ht!]
\begin{center}
\begin{tikzpicture}[scale=0.45,style=thick]
\def\vr{4pt}
\path (-10,3) coordinate (a'1);
\path (-10,1) coordinate (a''1);
\path (-12,2) coordinate (b'1);
\path (-6,2) coordinate (c'1);
\path (-8,3) coordinate (d'1);
\path (-8,1) coordinate (d''1);
\path (-9,5) coordinate (e1);
\path (-10,11) coordinate (a'2);
\path (-10,9) coordinate (a''2);
\path (-12,10) coordinate (b'2);
\path (-6,10) coordinate (c'2);
\path (-8,11) coordinate (d'2);
\path (-8,9) coordinate (d''2);
\path (-9,7) coordinate (e2);


\path (8,8) coordinate (e);
\path (7,11) coordinate (a);
\path (7,13) coordinate (aa);
\path (6,12) coordinate (b);
\path (10,12) coordinate (c);
\path (9,11) coordinate (d);
\path (9,13) coordinate (dd);

\path (8,4) coordinate (e');
\path (7,1) coordinate (a');
\path (7,-1) coordinate (aa');
\path (6,0) coordinate (b');
\path (10,0) coordinate (c');
\path (9,1) coordinate (d');
\path (9,-1) coordinate (dd');

\path (2,6) coordinate (e'');
\path (1,3) coordinate (a'');
\path (1,1) coordinate (aa'');
\path (0,2) coordinate (b'');
\path (4,2) coordinate (c'');
\path (3,3) coordinate (d'');
\path (3,1) coordinate (dd'');
\path (6,6) coordinate (x);


 \draw (e)--(x)--(e');
\draw (x)--(e'');

\draw (dd)--(a) -- (b); 
\draw (d)--(c);
\draw (a) -- (d)--(aa); 
\draw (b)--(aa)--(dd)--(c);
\draw (b)--(e)--(c);

\draw (dd')--(a') -- (b'); 
\draw (d')--(c');
\draw (a') -- (d')--(aa'); 
\draw (b')--(aa')--(dd')--(c');
\draw (b')--(e')--(c');

\draw (dd'')--(a'') -- (b''); 
\draw (d'')--(c'');
\draw (a'') -- (d'')--(aa''); 
\draw (b'')--(aa'')--(dd'')--(c'');
\draw (b'')--(e'')--(c'');

\draw (e2)--(e1);
\draw (d''2)--(a'2) -- (b'2); 
\draw (d'2)--(c'2);
\draw (a'2) -- (d'2)--(a''2); 
\draw (b'2)--(a''2)--(d''2)--(c'2);
\draw (b'2)--(e2)--(c'2);

\draw (d''1)--(a'1) -- (b'1); 
\draw (d'1)--(c'1);
\draw (a'1) -- (d'1)--(a''1); 
\draw (b'1)--(a''1)--(d''1)--(c'1);
\draw (b'1)--(e1)--(c'1);

\draw (a) [fill=white] circle (\vr);
\draw (aa) [fill=white] circle (\vr);
\draw (b) [fill=white] circle (\vr);
\draw (c) [fill=white] circle (\vr);
\draw (d) [fill=white] circle (\vr);
\draw (dd) [fill=white] circle (\vr);
\draw (e) [fill=white] circle (\vr);

\draw (a') [fill=white] circle (\vr);
\draw (aa') [fill=white] circle (\vr);
\draw (b') [fill=white] circle (\vr);
\draw (c') [fill=white] circle (\vr);
\draw (d') [fill=white] circle (\vr);
\draw (dd') [fill=white] circle (\vr);
\draw (e') [fill=white] circle (\vr);

\draw (a'') [fill=white] circle (\vr);
\draw (aa'') [fill=white] circle (\vr);
\draw (b'') [fill=white] circle (\vr);
\draw (c'') [fill=white] circle (\vr);
\draw (d'') [fill=white] circle (\vr);
\draw (dd'') [fill=white] circle (\vr);
\draw (e'') [fill=white] circle (\vr);
\draw (x) [fill=white] circle (\vr);

\draw (a'2) [fill=white] circle (\vr);
\draw (a''2) [fill=white] circle (\vr);
\draw (b'2) [fill=white] circle (\vr);
\draw (c'2) [fill=white] circle (\vr);
\draw (d'2) [fill=white] circle (\vr);
\draw (d''2) [fill=white] circle (\vr);
\draw (e2) [fill=white] circle (\vr);

\draw (a'1) [fill=white] circle (\vr);
\draw (a''1) [fill=white] circle (\vr);
\draw (b'1) [fill=white] circle (\vr);
\draw (c'1) [fill=white] circle (\vr);
\draw (d'1) [fill=white] circle (\vr);
\draw (d''1) [fill=white] circle (\vr);
\draw (e1) [fill=white] circle (\vr);

\end{tikzpicture}
\end{center}
\caption{Graphs $X_2$ and $X_3$.}
\label{fig:X2X3}
\end{figure}
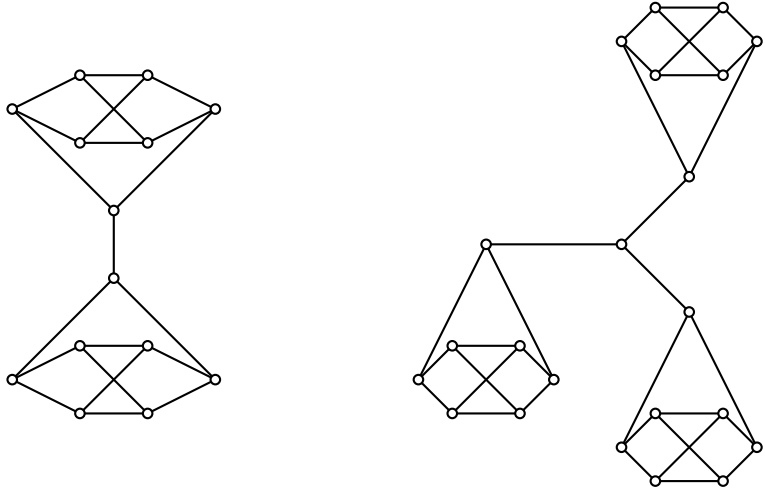

\begin{figure}[ht!]
\begin{center}
\begin{tikzpicture}[scale=0.45,style=thick]
\def\vr{4pt}
\path (-15,10) coordinate (xa);
\path (-17,9) coordinate (xb);
\path (-13,9) coordinate (xc);
\path (-15,8) coordinate (xd);
\path (-15,6) coordinate (xe);
\path (-15,4) coordinate (xe');
\path (-15.6,2) coordinate (xa');
\path (-15.6,0) coordinate (xa'');
\path (-17,1) coordinate (xb');
\path (-13,1) coordinate (xc');
\path (-14.4,2) coordinate (xd');
\path (-14.4,0) coordinate (xd'');

\draw (xa) -- (xb) -- (xd)--(xc)--(xa); 
\draw (xa) -- (xd);
\draw (xb)--(xe)--(xc); 
\draw (xd'')--(xa') -- (xb'); 
\draw (xd')--(xc');
\draw (xa') -- (xd')--(xa'');
\draw (xb')--(xe')--(xc'); 
\draw (xe)--(xe');
\draw (xb')--(xa'')--(xd'')--(xc');

\draw (xa) [fill=white] circle (\vr);
\draw (xb) [fill=white] circle (\vr);
\draw (xc) [fill=white] circle (\vr);
\draw (xd) [fill=white] circle (\vr);
\draw (xe) [fill=white] circle (\vr);
\draw (xa') [fill=white] circle (\vr);
\draw (xa'') [fill=white] circle (\vr);
\draw (xb') [fill=white] circle (\vr);
\draw (xc') [fill=white] circle (\vr);
\draw (xd') [fill=white] circle (\vr);
\draw (xd'') [fill=white] circle (\vr);
\draw (xe') [fill=white] circle (\vr);

\path (-6,10) coordinate (a2);
\path (-8,9) coordinate (b2);
\path (-4,9) coordinate (c2);
\path (-6,8) coordinate (d2);
\path (-6,6) coordinate (e2);
\path (-6,4) coordinate (e'2);
\path (-6,2) coordinate (a'2);
\path (-8,1) coordinate (b'2);
\path (-4,1) coordinate (c'2);
\path (-6,0) coordinate (d'2);
\path (8,12) coordinate (a);
\path (10,11) coordinate (b);
\path (6,11) coordinate (c);
\path (8,10) coordinate (d);
\path (8,8) coordinate (e);
\path (8,4) coordinate (e');
\path (8,2) coordinate (a');
\path (10,1) coordinate (b');
\path (6,1) coordinate (c');
\path (8,0) coordinate (d');

\path (2,6) coordinate (e'');
\path (2,4) coordinate (a'');
\path (4,3) coordinate (b'');
\path (0,3) coordinate (c'');
\path (2,2) coordinate (d'');

\path (6,6) coordinate (x);

\draw (a) -- (b) -- (d)--(c)--(a); 
\draw (a) -- (d);
\draw (b)--(e)--(c); 
\draw (a') -- (b') -- (d')--(c')--(a'); 
\draw (a') -- (d');
\draw (b')--(e')--(c');
\draw (a'') -- (b'') -- (d'')--(c'')--(a''); 
\draw (a'') -- (d'');
 \draw (e)--(x)--(e');
\draw (x)--(e'')--(c'');
\draw (e'')--(b'');

\draw (a2) -- (b2) -- (d2)--(c2)--(a2); 
\draw (a2) -- (d2);
\draw (b2)--(e2)--(c2); 
\draw (a'2) -- (b'2) -- (d'2)--(c'2)--(a'2); 
\draw (a'2) -- (d'2);
\draw (b'2)--(e'2)--(c'2); 
\draw (e2)--(e'2);

\draw (a) [fill=white] circle (\vr);
\draw (b) [fill=white] circle (\vr);
\draw (c) [fill=white] circle (\vr);
\draw (d) [fill=white] circle (\vr);
\draw (e) [fill=white] circle (\vr);
\draw (a') [fill=white] circle (\vr);
\draw (b') [fill=white] circle (\vr);
\draw (c') [fill=white] circle (\vr);
\draw (d') [fill=white] circle (\vr);
\draw (e') [fill=white] circle (\vr);
\draw (a'') [fill=white] circle (\vr);
\draw (b'') [fill=white] circle (\vr);
\draw (c'') [fill=white] circle (\vr);
\draw (d'') [fill=white] circle (\vr);
\draw (e'') [fill=white] circle (\vr);
\draw (x) [fill=white] circle (\vr);

\draw (a2) [fill=white] circle (\vr);
\draw (b2) [fill=white] circle (\vr);
\draw (c2) [fill=white] circle (\vr);
\draw (d2) [fill=white] circle (\vr);
\draw (e2) [fill=white] circle (\vr);
\draw (a'2) [fill=white] circle (\vr);
\draw (b'2) [fill=white] circle (\vr);
\draw (c'2) [fill=white] circle (\vr);
\draw (d'2) [fill=white] circle (\vr);
\draw (e'2) [fill=white] circle (\vr);

\end{tikzpicture}
\end{center}
\caption{Graphs $XY$, $Y_2$ and $Y_3$.}
\label{fig:Y2Y3}
\end{figure}
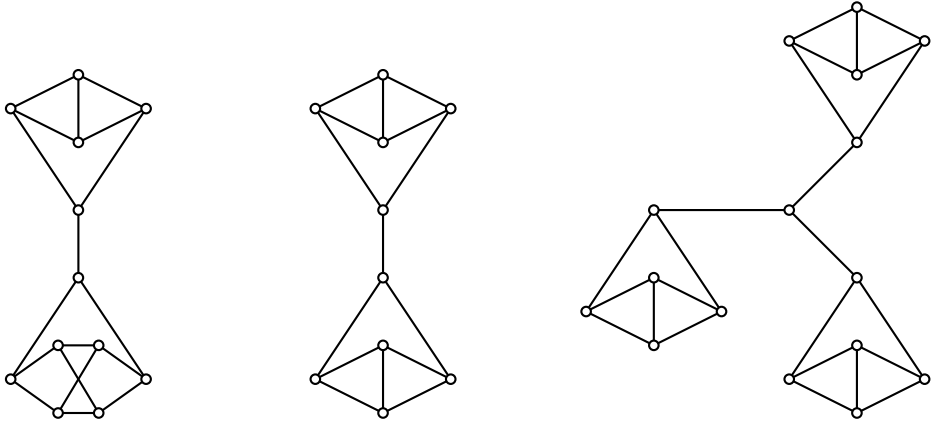

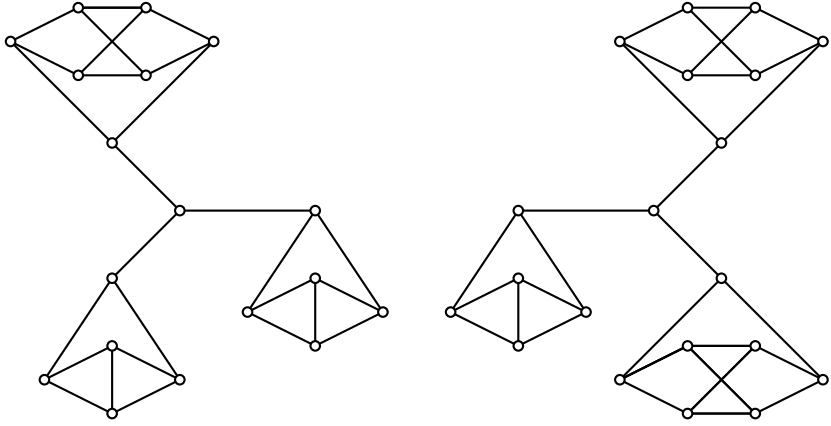
\begin{figure}[ht!]
\begin{center}
\begin{tikzpicture}[scale=0.45,style=thick]
\def\vr{4pt}

\path (-9,12) coordinate (1a);
\path (-9,10) coordinate (1aa);
\path (-11,11) coordinate (1b);
\path (-5,11) coordinate (1c);
\path (-7,12) coordinate (1d);
\path (-7,10) coordinate (1dd);
\path (-8,8) coordinate (1e);
\path (-8,4) coordinate (1e');
\path (-8,2) coordinate (1a');
\path (-10,1) coordinate (1b');
\path (-6,1) coordinate (1c');
\path (-8,0) coordinate (1d');

\path (-2,6) coordinate (1e'');
\path (-2,4) coordinate (1a'');
\path (-4,3) coordinate (1b'');
\path (-0,3) coordinate (1c'');
\path (-2,2) coordinate (1d'');

\path (-6,6) coordinate (1x);

\path (9,12) coordinate (aa);
\path (9,10) coordinate (a);
\path (13,11) coordinate (b);
\path (7,11) coordinate (c);
\path (11,12) coordinate (d);
\path (11,10) coordinate (dd);

\path (10,8) coordinate (e);
\path (10,4) coordinate (e');

\path (9,2) coordinate (a');
\path (9,0) coordinate (a'a);
\path (13,1) coordinate (b');
\path (7,1) coordinate (c');
\path (11,0) coordinate (d'd);
\path (11,2) coordinate (d');

\path (4,6) coordinate (e'');
\path (4,4) coordinate (a'');
\path (6,3) coordinate (b'');
\path (2,3) coordinate (c'');
\path (4,2) coordinate (d'');

\path (8,6) coordinate (x);

\draw (1aa)--(1b) -- (1a)--(1d)--(1c); 
\draw (1dd)--(1a) -- (1d)--(1aa)--(1dd);
\draw (1b)--(1e)--(1c)--(1dd); 
\draw (1a') -- (1b') -- (1d')--(1c')--(1a'); 
\draw (1a') -- (1d');
\draw (1b')--(1e')--(1c');
\draw (1a'') -- (1b'') -- (1d'')--(1c'')--(1a''); 
\draw (1a'') -- (1d'');
\draw (1e)--(1x)--(1e');
\draw (1x)--(1e'')--(1c'');
\draw (1e'')--(1b'');

\draw (dd)--(b) -- (d);
\draw (a)--(c)--(aa); 
\draw (a) -- (d)--(aa)--(dd)-- (a);
\draw (b)--(e)--(c); 
\draw (d'd)--(b') -- (d');
\draw (a')--(c')--(a'a); 
\draw (a') -- (d')--(a'a)--(d'd)-- (a');
\draw (a') -- (c');
\draw (c')--(a'); 
\draw (a') -- (d')--(a'a)--(d'd)-- (a');
\draw (b')--(e')--(c');

\draw (a'') -- (b'') -- (d'')--(c'')--(a''); 
\draw (a'') -- (d'');
 \draw (e)--(x)--(e');
\draw (x)--(e'')--(c'');
\draw (e'')--(b'');

\draw (a) [fill=white] circle (\vr);
\draw (b) [fill=white] circle (\vr);
\draw (c) [fill=white] circle (\vr);
\draw (d) [fill=white] circle (\vr);
\draw (aa) [fill=white] circle (\vr);
\draw (a'a) [fill=white] circle (\vr);
\draw (dd) [fill=white] circle (\vr);
\draw (d'd) [fill=white] circle (\vr);
\draw (e) [fill=white] circle (\vr);
\draw (a') [fill=white] circle (\vr);
\draw (b') [fill=white] circle (\vr);
\draw (c') [fill=white] circle (\vr);
\draw (d') [fill=white] circle (\vr);
\draw (e') [fill=white] circle (\vr);
\draw (a'') [fill=white] circle (\vr);
\draw (b'') [fill=white] circle (\vr);
\draw (c'') [fill=white] circle (\vr);
\draw (d'') [fill=white] circle (\vr);
\draw (e'') [fill=white] circle (\vr);
\draw (x) [fill=white] circle (\vr);

\draw (1a) [fill=white] circle (\vr);
\draw (1aa) [fill=white] circle (\vr);
\draw (1b) [fill=white] circle (\vr);
\draw (1c) [fill=white] circle (\vr);
\draw (1d) [fill=white] circle (\vr);
\draw (1dd) [fill=white] circle (\vr);
\draw (1e) [fill=white] circle (\vr);
\draw (1a') [fill=white] circle (\vr);
\draw (1b') [fill=white] circle (\vr);
\draw (1c') [fill=white] circle (\vr);
\draw (1d') [fill=white] circle (\vr);
\draw (1e') [fill=white] circle (\vr);
\draw (1a'') [fill=white] circle (\vr);
\draw (1b'') [fill=white] circle (\vr);
\draw (1c'') [fill=white] circle (\vr);
\draw (1d'') [fill=white] circle (\vr);
\draw (1e'') [fill=white] circle (\vr);
\draw (1x) [fill=white] circle (\vr);

\end{tikzpicture}
\end{center}
\caption{Graphs $XY_2$ and $X_2Y$.}
\label{fig:XY2}
\end{figure}

Taking two copies of the graph $X$ and connecting with an edge the vertices of degree $2$ creates a cubic graph, which we denote by $X_2$ (it is obtained from the tree $K_2$ by using the above identification of leaves with the vertices $\alpha_X$ of their own copy of $X$).
Figure~\ref{fig:X2X3} shows the graphs $X_2$ and $X_3$, where the latter is obtained by the described operation from the tree $K_{1,3}$. In a similar way we construct the graphs $Y_2$ and $Y_3$, see Figure~\ref{fig:Y2Y3}. On the left of this figure the graph $XY$ is depicted, which is obtained from a copy of $X$ and a copy of $Y$ by adding an edge between $\alpha_X$ and $\alpha_Y$.  
The graphs obtained from $K_{1,3}$ by attaching to the three leaves two copies of $X$ and a copy of $Y$ or two copies of $Y$ and a copy of $X$ are denoted by $X_2Y$ and $XY_2$, respectively. They are depicted in Figure~\ref{fig:XY2}. 

We follow with a characterization of the graphs in family $\cal M$ that attain the bound in Theorem~\ref{thm:Z-cubic}. The subfamily ${\cal M}'$ of $\cal M$ consists of the graphs  $X_2, X_3,Y_2,Y_3,X_2Y,XY_2$ and $XY$.

\begin{prop}
\label{prp:familyM}
If $G\in \cal M$, then $\grz(G)=\frac{n(G)}{2}$ if and only if $G\in {\cal M}'$. 
\end{prop}
\proof It is easy to see that $\grz(X_2)=7$, $\grz(X_3)=11$, $\grz(Y_2)=5$, $\grz(Y_3)=8$, $\grz(XY)=6$, $\grz(XY_2)=9$ and $\grz(X_2Y)=10$. This shows that $G\in {\cal M}'$ implies $\grz(G)=\frac{n(G)}{2}$.

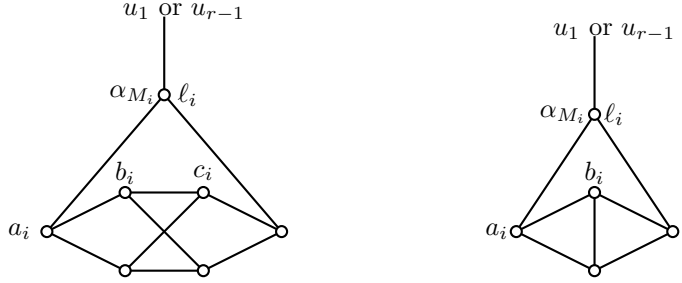
\begin{figure}[ht!]
\begin{center}
\begin{tikzpicture}[scale=0.52,style=thick]
\def\vr{4pt}

\path (6,6) coordinate (h);
\path (-5,6.5) coordinate (g);

\path (6,4) coordinate (e'2);
\path (6,2) coordinate (a'2);
\path (8,1) coordinate (b'2);
\path (4,1) coordinate (c'2);
\path (6,0) coordinate (d'2);
\path (-6,2) coordinate (a');
\path (-6,0) coordinate (a'');
\path (-8,1) coordinate (b');
\path (-2,1) coordinate (c');
\path (-4,2) coordinate (d');
\path (-4,0) coordinate (d'');
\path (-5,4.5) coordinate (e);

\draw (5.2,4) node [text=black]{$\alpha_{M_i}$}; 
\draw (6.5,4) node [text=black]{$\ell_i$}; 
\draw (-5.8,4.5) node [text=black]{$\alpha_{M_i}$};
\draw (-4.4,4.5) node [text=black]{$\ell_i$}; 
\draw (-4.5,6.6) node [text=black]{$u_1\textrm{ or }u_{r-1}$}; 
\draw (6.5,6.1) node [text=black]{$u_1\textrm{ or }u_{r-1}$}; 

\draw (-8.7,1) node [text=black]{$a_i$}; 
\draw (-6,2.5) node [text=black]{$b_i$}; 
\draw (-4,2.5) node [text=black]{$c_i$}; 

\draw (3.5,1) node [text=black]{$a_i$}; 
\draw (6,2.5) node [text=black]{$b_i$};

\draw (a'2) -- (b'2) -- (d'2)--(c'2)--(a'2); 
\draw (a'2) -- (d'2);
\draw (b'2)--(e'2)--(c'2); 
\draw (e'2)--(h);
\draw (g)--(e);


\draw (d'')--(a') -- (b'); 
\draw (d')--(c');
\draw (a') -- (d')--(a''); 
\draw (b')--(a'')--(d'')--(c');
\draw (b')--(e)--(c');

\draw (a'2) [fill=white] circle (\vr);
\draw (b'2) [fill=white] circle (\vr);
\draw (c'2) [fill=white] circle (\vr);
\draw (d'2) [fill=white] circle (\vr);
\draw (e'2) [fill=white] circle (\vr);


\draw (a') [fill=white] circle (\vr);
\draw (a'') [fill=white] circle (\vr);
\draw (b') [fill=white] circle (\vr);
\draw (c') [fill=white] circle (\vr);
\draw (d') [fill=white] circle (\vr);
\draw (d'') [fill=white] circle (\vr);
\draw (e) [fill=white] circle (\vr);

\end{tikzpicture}
\end{center}
\caption{Graph $M_i$ in the proof of Proposition~\ref{prp:familyM} 
as a copy of $X$ and a copy $Y$.}
\label{fig:XandYasM}
\end{figure}

For the other direction, we will construct a Z-sequence in a graph $G\in  {\cal M}\setminus {\cal M}'$ with length greater than $n(G)/2$. Note that $G$ can be obtained from a tree $T$ whose non-leaf vertices have degree $3$ by identifying each leaf $\ell$ of $T$ with the vertex $x_M$ or $y_M$ that corresponds to $\ell$. In addition, since $G$ is not in ${\cal M}'$, the tree $T$ has more than four vertices and contains at least two leaves, say $\ell_1$ and $\ell_2$, that are at distance greater than $2$. Consider the two copies of graphs in $\{X,Y\}$ that are associated to the leaves $\ell_1$ and $\ell_2$, and denote them by $M_1$ and $M_2$, respectively; similarly, the identified vertices of $M_1$ and $M_2$ are denoted by $\alpha_{M_1}$ and $\alpha_{M_2}$, respectively. Note that each of the graphs $M_i$ can be isomorphic either to $X$ or $Y$. In Figure~\ref{fig:XandYasM}, a graph $M_i$ is depicted, on the left as a copy of $X$ and on the right as a copy of $Y$. Consider also the notation of vertices on the same figure.

Let $\ell_1=u_0,u_1\ldots,u_r=\ell_2$ be the shortest path in $T$ (and also in $G$) between $\alpha_{M_1}=\ell_1$ and $\alpha_{M_2}=\ell_2$. Note that $r\ge 3$. Now, depending on which of the graphs are $M_1$ and $M_2$, consider a sequence $S$, which starts in the subgraphs $M_1$ and $M_2$ as follows. If both $M_i$ are isomorphic to $X$, then $S$ starts with $(a_1,b_1,c_1,\alpha_{M_1}, a_2,b_2,c_2,\alpha_{M_2},u_1,u_2,\ldots,u_{r-1})$, and note that $a_1$ and $a_2$ footprint four vertices, $c_1, c_2,\alpha_{M_1},\alpha_{M_2},u_{r-2}$ and $u_{r-1}$ footprint only one vertex, while all other vertices footprint two vertices. In the same way as in the proof of Lemma~\ref{lem:ZsequenceVec} we note that the Z-sequence $S$ can be extended to a dominating Z-sequence of length $t$ such that each the remaining vertices of the sequence footprints at most two vertices. We infer that the vertices of the sequence footprint at most $4+4+6+2(t-8)$ vertices, which implies that $n\le 4+4+6+2(t-8)$, and so $$\grz(G)\ge t\ge \frac{n+2}{2}>\frac{n(G)}{2}.$$ 

The case when both $M_i$ are isomorphic to $Y$, or one is isomorphic to $X$ and the other to $Y$, can be proved in a similar way. For instance, in the former case, the appropriate sequence $S$ starts with $(a_1,b_1,\alpha_{M_1}, a_2,b_2,\alpha_{M_2},u_1,u_2,\ldots,u_{r-1})$, and again there are six vertices that footprint only one vertex. The remainder of the proof is essentially the same as above.\qed

\bigskip 

The following lemma indicates a special role that is played by graphs $X$ and $Y$ in the study of Z-Grundy domination in cubic graphs. It will be used several times in the proof of the main theorem. 

\begin{lemma}
\label{lem:Zlata}
Let $G$ be a connected cubic graph, $uv$ a cut edge in $G$, and let $u$ lie on a cycle. Denote by $U$ the component of $G-uv$ that contains $u$.  
\begin{enumerate}[(i)]
\item If for every dominating Z-sequence $S$ in which no vertex of $\widehat{S}\cap V(U)$ footprints more than two vertices, every vertex of  $\widehat{S}\cap V(U)$ footprints exactly two vertices, then $U$ is isomorphic to $X$ or $Y$. 
\item If $U$ is isomorphic to $X$ or $Y$, then there exists a dominating Z-sequence $S$ that starts with a vertex in $U$, two vertices of $\widehat{S}\cap V(U)$ footprint only one vertex, and all other vertices of $\widehat{S}\cap V(U)$, except the first vertex of $S$, footprint at most two vertices.
\end{enumerate}
\end{lemma}

\proof (i) 
Let $C$ be a cycle that contains $u$ whose order $p$ is as small as possible. In particular, $C$ is an induced cycle. Let $C:u,c_2,\ldots,c_p,u$. Assume first that $p=3$. That is, $u, c_2$ and $c_3$ form a triangle, and each of $c_2$ and $c_3$ has a neighbor outside $C$. Then, if the sequence $S$ starts with $(v,u,c_2)$, vertex $c_{2}$ footprints exactly one vertex. Clearly, in the same way as in the proof of Lemma~\ref{lem:ZsequenceVec}, $S$ can be extended to a dominating Z-sequence in which every vertex of $C$ footprints at most two vertices, which is a contradiction.

Let $p=4$. If $c_3$ has a neighbor outside $C$, which is not adjacent to $c_2$, then starting the sequence $S$ with $(v,u,c_2,c_3)$  vertex $c_{3}$ footprints exactly one vertex, yielding again a contradiction. We may thus assume that $c_3$ and $c_2$ have a common neighbor $w$, and by symmetry (replacing the roles of $c_4$ and $c_2$), vertex $c_4$ is also adjacent to $w$.  In that case, vertices $u,c_2,c_3,c_4$ and $w$ induce a subgraph isomorphic to $Y$. 

Now, let $p\ge 5$. Suppose that the neighbor $w$ of $c_{p-1}$, which does not lie on $C$, is adjacent only to $c_{p-1}$ among all vertices of $C$. Then in the sequence $(v, u,c_2,\ldots,c_{p-1})$, as the starting part of $S$, $c_{p-1}$ footprints exactly one vertex, which yields the same contradiction as earlier. On the other hand, assume that $w$ is a common neighbor of another vertex from $C$. By the choice of $C$ as a shortest possible cycle containing $u$, we infer that $w$ can only be adjacent to $c_p$, $c_{p-2}$ or $c_{p-3}$. If $w$ is adjacent to $c_p$, then by the minimality of $C$ it cannot be adjacent also to $c_{p-3}$; also $w$ cannot be adjacent to $c_{p-2}$, since then in the sequence $(v,u,c_p,c_{p-1})$, the vertex $c_{p-1}$ footprints only $c_{p-2}$. On the other hand, if $w$ is not adjacent to $c_{p-2}$, then in the sequence $(v,u,c_p,w)$ the vertex $w$ footprints only one vertex. We thus infer that $wc_p\notin E(G)$. By $w'$ we denote the neighbor of $c_{p-2}$, which does not lie on $C$. 

Assume that $wc_{p-3}\in E(G)$. If  $w'c_p\notin E(G)$, then in $(v,u,c_2,\ldots,c_{p-2},c_p)$ vertex $c_p$ footprints exactly one vertex. Hence, let $w'c_p\in E(G)$.  (Note that then $w \neq w'$, since $C$ is a shortest cycle containing $u$.)  If $ww' \notin E(G)$, then in the sequence $(v,u,c_2,\ldots,c_{p-2},w')$ the vertex $w'$ footprints exactly one vertex. Now, let $ww'\in E(G)$. If $p \geq 6$, then in the sequence $(v,u,c_p,c_{p-1},c_{p-2})$, vertex $c_{p-2}$ footprints exactly one vertex. Finally, if $p=5$, then $C$ is isomorphic to $X$. 

(ii) The proof is straightforward: by using the notation from Figure~\ref{fig:XandYasM}, the corresponding sequence for graph $X$ 
starts with $(a_i,b_i,c_i,\alpha_{M_i})$, and corresponding sequence for graph $Y$ starts with $(a_i,b_i,\alpha_{M_i})$. (Note that the role of vertex $u$ is played by $\alpha_{M_i}$.)
\qed

\begin{figure}[ht!]
\begin{center}
\begin{tikzpicture}[scale=0.52,style=thick]
\def\vr{4pt}

\path (-1,1.5) coordinate (2x);
\path (-3,0) coordinate (2y);
\path (-5,1.5) coordinate (2z);
\path (-1,5) coordinate (2u);
\path (-3,5) coordinate (2v);
\path (-2,7) coordinate (2v');
\path (-4,7) coordinate (2v'');
\path (-5,5) coordinate (2w);


\draw (2x) -- (2y) -- (2z) -- (2x); 
\draw (2w) -- (2v') -- (2v) -- (2v'') -- (2u);
\draw (2v') -- (2u); 
\draw (2v'') -- (2w); 
\draw (2x) -- (2u);
\draw (2y) -- (2v);  
\draw (2z) -- (2w);


\draw (2x) [fill=white] circle (\vr);
\draw (2y) [fill=white] circle (\vr);
\draw (2z) [fill=white] circle (\vr);
\draw (2u) [fill=white] circle (\vr);
\draw (2v) [fill=white] circle (\vr);
\draw (2v') [fill=white] circle (\vr);
\draw (2v'') [fill=white] circle (\vr);
\draw (2w) [fill=white] circle (\vr);

\path (7,1.5) coordinate (x);
\path (9,0) coordinate (y);
\path (11,1.5) coordinate (z);
\path (7,6.5) coordinate (u);
\path (9,5) coordinate (v);
\path (11,6.5) coordinate (w);


\draw (x) -- (y) -- (z) -- (x); 
\draw (u) -- (v) -- (w) -- (u);
\draw (x) -- (u);
\draw (y) -- (v);  
\draw (z) -- (w);


\draw (x) [fill=white] circle (\vr);
\draw (y) [fill=white] circle (\vr);
\draw (z) [fill=white] circle (\vr);
\draw (u) [fill=white] circle (\vr);
\draw (v) [fill=white] circle (\vr);
\draw (w) [fill=white] circle (\vr);

\end{tikzpicture}
\end{center}
\caption{Graphs $TK$ and $K_3\Box K_2$.}
\label{fig:hamming}
\end{figure}
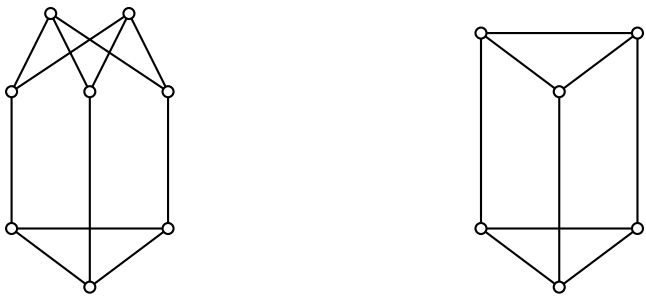

\bigskip
We next present several special (cubic) graphs that appear in the characterization. We start with the cubic graph on $8$ vertices, which we denote by $TK$, and is depicted on the left side of Figure~\ref{fig:hamming}. On the right 
side of the same figure one can see the Hamming graph $K_3\cp K_2$.

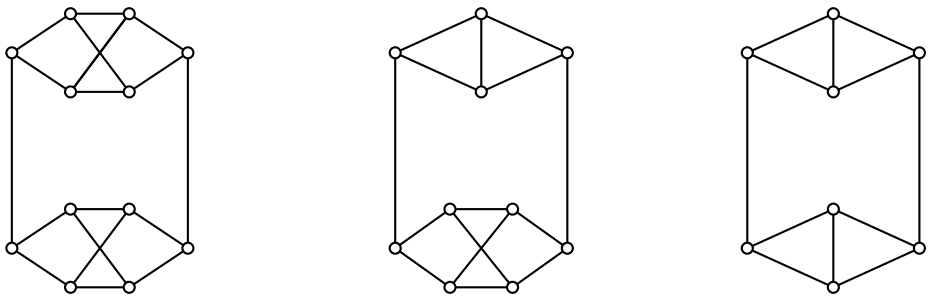
\begin{figure}[ht!]
\begin{center}
\begin{tikzpicture}[scale=0.52,style=thick]
\def\vr{4pt}
\path (-9,7) coordinate (c);
\path (-9,5) coordinate (a);
\path (-10.5,7) coordinate (d);
\path (-10.5,5) coordinate (b);

\path (-10.5,2) coordinate (a');
\path (-10.5,0) coordinate (a'');

\path (-12,6) coordinate (e);
\path (-7.5,6) coordinate (e');

\path (-12,1) coordinate (b');
\path (-7.5,1) coordinate (c');
\path (-9,2) coordinate (d');
\path (-9,0) coordinate (d'');

\draw (a) -- (b); 
\draw (a) -- (d); 
\draw (a) -- (e');
\draw (b) -- (c);
\draw (b) -- (e);
\draw (c)--(e'); 
\draw (c) -- (d);
\draw (c) -- (b);
\draw (d) -- (e); 
\draw (d'')--(a') -- (b'); 
\draw (d')--(c');
\draw (e')--(c');
\draw (e)--(b');
\draw (a') -- (d')--(a''); 
\draw (b')--(a'')--(d'')--(c');

\draw (a) [fill=white] circle (\vr);
\draw (b) [fill=white] circle (\vr);
\draw (c) [fill=white] circle (\vr);
\draw (d) [fill=white] circle (\vr);
\draw (e) [fill=white] circle (\vr);
\draw (a') [fill=white] circle (\vr);
\draw (a'') [fill=white] circle (\vr);
\draw (b') [fill=white] circle (\vr);
\draw (c') [fill=white] circle (\vr);
\draw (d') [fill=white] circle (\vr);
\draw (d'') [fill=white] circle (\vr);
\draw (e') [fill=white] circle (\vr);


\path (0,7) coordinate (k);
\path (-2.2,6) coordinate (l);
\path (2.2,6) coordinate (m);
\path (0,5) coordinate (n);

\path (-0.8,2) coordinate (k');
\path (-0.8,0) coordinate (k'');
\path (-2.2,1) coordinate (l');
\path (2.2,1) coordinate (m');
\path (0.8,2) coordinate (n');
\path (0.8,0) coordinate (n'');

\draw (k) -- (l) -- (n)--(m)--(k); 
\draw (k) -- (n); 
\draw (k') -- (l')--(k''); 
\draw (n')--(m')--(n''); 
\draw (l) -- (l'); 
\draw (m) -- (m'); 
\draw (n'')--(k')--(n');
\draw (n'')--(k'')--(n');


\draw (k) [fill=white] circle (\vr);
\draw (l) [fill=white] circle (\vr);
\draw (m) [fill=white] circle (\vr);
\draw (n) [fill=white] circle (\vr);
\draw (k') [fill=white] circle (\vr);
\draw (k'') [fill=white] circle (\vr);
\draw (l') [fill=white] circle (\vr);
\draw (m') [fill=white] circle (\vr);
\draw (n') [fill=white] circle (\vr);
\draw (n'') [fill=white] circle (\vr);


\path (9,7) coordinate (2k);
\path (6.8,6) coordinate (2l);
\path (11.2,6) coordinate (2m);
\path (9,5) coordinate (2n);

\path (9,2) coordinate (2k');
\path (6.8,1) coordinate (2l');
\path (11.2,1) coordinate (2m');
\path (9,0) coordinate (2n');

\draw (2k) -- (2l) -- (2n)--(2m)--(2k); 
\draw (2k) -- (2n); 
\draw (2k') -- (2l') -- (2n')--(2m')--(2k'); 
\draw (2k') -- (2n'); 
\draw (2l) -- (2l'); 
\draw (2m) -- (2m');

\draw (2k) [fill=white] circle (\vr);
\draw (2l) [fill=white] circle (\vr);
\draw (2m) [fill=white] circle (\vr);
\draw (2n) [fill=white] circle (\vr);
\draw (2k') [fill=white] circle (\vr);
\draw (2l') [fill=white] circle (\vr);
\draw (2m') [fill=white] circle (\vr);
\draw (2n') [fill=white] circle (\vr);

\end{tikzpicture}
\end{center}
\caption{Graphs necklaces $N_{XX}$, $N_{XY}$, and $N_{YY}$.}
\label{fig:necklace}
\end{figure}

The {\em diamond} is the graph $K_4-e$ obtained from the complete graph on $4$ vertices by deleting an edge. 
We next define three cubic graphs, which we will refer to as {\em necklaces}. We denote them by $N_{XX}, N_{XY}$ and $N_{YY}$, see Figure~\ref{fig:necklace}. In particular, the cubic graph obtained from the disjoint union of two copies of the diamond by adding two edges is known as the 2-necklace, it is depicted on the right side of Figure~\ref{fig:necklace}, and we will denote it by $N_{YY}$.

\bigskip

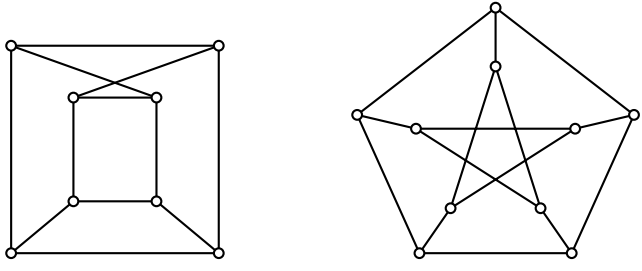
\begin{figure}[ht!]
\begin{center}
\begin{tikzpicture}[scale=0.46,style=thick]
\def\vr{4pt}
\path (-2,0) coordinate (a);
\path (-8,0) coordinate (b);
\path (-2,6) coordinate (c);
\path (-8,6) coordinate (d);
\path (-3.8,1.5) coordinate (e);
\path (-6.2,1.5) coordinate (f);
\path (-3.8,4.5) coordinate (g);
\path (-6.2,4.5) coordinate (h);

\path (3.8,0) coordinate (x);
\path (8.2,0) coordinate (y);
\path (2,4) coordinate (z);
\path (10,4) coordinate (u);
\path (6,7.1) coordinate (v);
\path (4.7,1.3) coordinate (x');
\path (7.3,1.3) coordinate (y');
\path (3.7,3.6) coordinate (z');
\path (8.3,3.6) coordinate (u');
\path (6,5.4) coordinate (v');

\draw (a) -- (b) -- (d)--(c)--(a); 
\draw (e) -- (g) --(h)--(f)--(e);
\draw (a) -- (e); 
\draw (b) -- (f); 
\draw (c) -- (h); 
\draw (d) -- (g); 

\draw (x) -- (y) -- (u)--(v)--(z) --(x); 
\draw (x) -- (x');
\draw (y) -- (y');
\draw (z) -- (z');
\draw (u) -- (u');
\draw (v) -- (v');
\draw (x') -- (u') -- (z')--(y')--(v') --(x'); 

\draw (a) [fill=white] circle (\vr);
\draw (b) [fill=white] circle (\vr);
\draw (c) [fill=white] circle (\vr);
\draw (d) [fill=white] circle (\vr);
\draw (e) [fill=white] circle (\vr);
\draw (f) [fill=white] circle (\vr);
\draw (g) [fill=white] circle (\vr);
\draw (h) [fill=white] circle (\vr);
\draw (x) [fill=white] circle (\vr);
\draw (y) [fill=white] circle (\vr);
\draw (z) [fill=white] circle (\vr);
\draw (u) [fill=white] circle (\vr);
\draw (v) [fill=white] circle (\vr);
\draw (x') [fill=white] circle (\vr);
\draw (y') [fill=white] circle (\vr);
\draw (z') [fill=white] circle (\vr);
\draw (u') [fill=white] circle (\vr);
\draw (v') [fill=white] circle (\vr);

\end{tikzpicture}
\end{center}
\caption{Twisted cube $TQ_3$ and the Petersen graph.}
\label{fig:non-bip}
\end{figure}

We follow with two non-bipartite triangle-free graphs.
The {\em twisted cube} $TQ_3$ is the graph obtained from the disjoint union of two $4$-cycles, say $C:a,b,c,d,a$ and $C':a',b',c',d',a'$ by adding the edges $aa'$, $bb'$, $cd'$ and $dc'$. Clearly, $TQ_3$ is a non-bipartite cubic graph of order $8$; see Figure~\ref{fig:non-bip}. On the right side of the same figure the Petersen graph is shown. 

\begin{figure}[ht!]
\begin{center}
\begin{tikzpicture}[scale=0.46,style=thick]
\def\vr{4pt}
\path (-2,0) coordinate (a);
\path (-8,0) coordinate (b);
\path (-2,6) coordinate (c);
\path (-8,6) coordinate (d);
\path (-3.8,1.5) coordinate (e);
\path (-6.2,1.5) coordinate (f);
\path (-3.8,4.5) coordinate (g);
\path (-6.2,4.5) coordinate (h);

\path (4,0.5) coordinate (x);
\path (7,0.5) coordinate (y);
\path (10,0.5) coordinate (z);
\path (4,5.5) coordinate (u);
\path (7,5.5) coordinate (v);
\path (10,5.5) coordinate (w);

\draw (a) -- (b) -- (d)--(c)--(a); 
\draw (e) -- (g) --(h)--(f)--(e);
\draw (a) -- (e); 
\draw (b) -- (f); 
\draw (c) -- (g); 
\draw (d) -- (h); 

\draw (x) -- (u) -- (y) -- (v) -- (z)--(w) --(x) --(v); 
\draw (y) -- (w);
\draw (z) -- (u);  

\draw (a) [fill=white] circle (\vr);
\draw (b) [fill=white] circle (\vr);
\draw (c) [fill=white] circle (\vr);
\draw (d) [fill=white] circle (\vr);
\draw (e) [fill=white] circle (\vr);
\draw (f) [fill=white] circle (\vr);
\draw (g) [fill=white] circle (\vr);
\draw (h) [fill=white] circle (\vr);
\draw (x) [fill=white] circle (\vr);
\draw (y) [fill=white] circle (\vr);
\draw (z) [fill=white] circle (\vr);
\draw (u) [fill=white] circle (\vr);
\draw (v) [fill=white] circle (\vr);
\draw (w) [fill=white] circle (\vr);

\end{tikzpicture}
\end{center}
\caption{Graphs $Q_3$ and $K_{3,3}$.}
\label{fig:bip}
\end{figure}
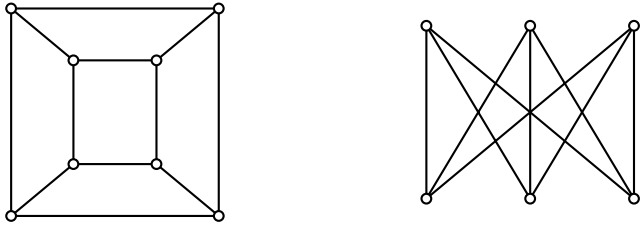

Finally, two bipartite graphs relevant for this section are the $3$-cube $Q_3$, and the complete bipartite graph $K_{3,3}$, shown in Figure~\ref{fig:bip}.

We are ready to formulate the characterization of the connected cubic graphs with Z-Grundy domination number half their order.

\begin{thm}
\label{thm:cubic}
A connected, cubic, graph $G$ of order $n$ has $\gr^Z(G)=\frac{n}{2}$ if and only if $G\in {\cal M}'$ or $G$ is one of the graphs $N_{XX}, N_{XY}, N_{YY}$, $K_3\Box K_2$, $TK$, $Q_3$, $TQ_3$, or the Petersen graph.
\end{thm}

\proof Clearly, by using Proposition~\ref{prp:familyM} for the class ${\cal M}'$, and noting that $\gr^Z(N_{XX})=6$, $\gr^Z(N_{XY})=5$, $\gr^Z(N_{YY})=4$, $\gr^Z(K_3\Box K_2)=3$, $\gr^Z(TK)=4$, $\gr^Z(Q_3)=4$, $\gr^Z(TQ_3)=4$ and $\gr^Z(P)=5$, where $P$ is Petersen graph, we see that the graphs that appear in the statement of the theorem attain the bound $\frac{n}{2}$ for their Z-Grundy domination numbers.

For the reverse direction, we distinguish four cases. The first case is that $G$ contains $Y$ as a subgraph, the second is that it does not contain $Y$, but contains a diamond as a subgraph (clearly, a diamond has to be an induced subgraph, since $G$ is cubic and $G\ne K_4$), and the third case is that $G$ has no diamond, but has a triangle. The final case is that $G$ is triangle-free. In the first two cases, let $G$ have an induced diamond with vertices $a_1,a_2,b$ and $c$, where $a_1$ and $a_2$ are not adjacent. 

\medskip

{\bf Case 1.} $G$ has $Y$ as a subgraph.

Let $a_1$ and $a_2$ have a common neighbor $\alpha_Y$ different from $b$ and $c$, hence the subset $\{a_1,a_2,b,c,\alpha_Y\}$ induces $Y$ as a subgraph.  Let a sequence start with $S'=(a_1,c,\alpha_Y)$. Note that $c$ and $\alpha_Y$ footprint only one vertex. By Lemma~\ref{lem:ZsequenceVec}, it suffices 
to extend $S'$ to a Z-sequence $S$ such that it contains at least one more vertex that footprints only one vertex and all other vertices of $S$ (except $a_1$) footprint at most two vertices of $G$.

Let $u' \notin \{a_1,a_2\}$ be the neighbor of $\alpha_Y$. It is possible that $u'$ lies on a cycle. Otherwise, $u'$ has two neighbors, $z$ and $w$, distinct from $\alpha_Y$, and the only path between $z$ and $w$ is the one passing $u'$. Since $G$ is cubic and finite, there exists a vertex $a$ in $G-u'z$ (and also in $G-u'w$) such that $a$ lies on a cycle, but all internal vertices on the path between $\alpha_Y$ and $a$ do not lie on a cycle. Denote by $A$ the set of all such vertices $a$ in $G$ that lie on a cycle and there is a path between $\alpha_Y$ and $a$ passing through $u'$ such that all internal vertices of the path do not lie on a cycle (note that if already $u'$ lies on a cycle, then $A=\{u'\}$). 

Take any $a\in A$, and let $a'$ be the vertex on the path between $\alpha_Y$ and $a$, which is adjacent to $a$. Note that $aa'$ is a cut-edge in $G$, and since $a$ lies on a cycle, we can apply Lemma~\ref{lem:Zlata}(i). We infer that unless the component $U$ of $G-aa'$ that contains $a$ is isomorphic to $X$ or $Y$, there is a dominating Z-sequence $S$ to which we can extend $S'$ such that every vertex of $\widehat{S}\cap V(U)$ footprints at most two vertices and at least one vertex of $\widehat{S}\cap V(U)$ footprints exactly one vertex.
Applying Lemma~\ref{lem:ZsequenceVec}, and considering the assumption $\grz(G)=\frac{n}{2}$, we infer that $G$ belongs to $\cal M$, and by Proposition~\ref{prp:familyM} we get $G\in {\cal M}'$.

\medskip
   
{\bf Case 2.} $G$ has a diamond, but has no $Y$ as a subgraph.

Let $G$ have a diamond with vertices $a_1,a_2,b$ and $c$, let $a_i$ be adjacent also to $v_i$, for $i\in\{1,2\}$, and $v_1\ne v_2$. If $v_1v_2\in E(G)$, then $(a_1,c,a_2,v_2)$ is a Z-sequence in which each of $c,a_2$ and $v_2$ footprints only one vertex. Thus, assume that $v_1v_2\notin E(G)$. Suppose that $d(v_1,v_2)=2$. First, let $v_1$ and $v_2$ have two common neighbors $x_1$ and $x_2$. If $x_1x_2\in E(G)$, then $G$ is isomorphic to the necklace $N_{YY}$ and we are done. Otherwise, if $x_1$ and $x_2$ are open twins with $N(x_i)=\{v_1,v_2,y\}$, then $(a_1,c,a_2,v_2,x_1,y)$ is a Z-sequence in which each of the vertices $c,a_2$ and $y$ footprints only one vertex, which again implies $\gr^Z(G)> \frac{n}{2}$. Now, if $x_1$ and $x_2$ are not open twins, then $(a_1,c,a_2,v_2,x_1)$ is a Z-sequence with $c,a_2$ and $x_1$ being such vertices that footprint only one vertex, and by using Lemma~\ref{lem:ZsequenceVec} we get a contradition to the assumption $\grz(G)=\frac{n}{2}$. If $v_1$ and $v_2$ have only one common neighbor, say $x$, then $(a_1,c,a_2,v_2,v_1)$ is a Z-sequence, in which each of $c,a_2$ and $v_1$ footprints only one vertex. 
We may thus assume that $d(v_1,v_2)\ge 3$.  

Suppose that there is a path 
between $v_1$ and $v_2$ that does not pass $a_1$ (and $a_2$).
Let $P:v_1=z_0,z_1,\ldots,z_r=v_2$, where $r\ge 3$, be a shortest such path.  Let $u_i\in N(v_i)\setminus\{a_i\}$ be the neighbor of $v_i$ that is not in $P$. Suppose that $u_i$, $i\in [2]$, is not dominated by a vertex from $V(P)\setminus\{v_i)$. Without loss of generality we may assume that $u_2$ is not adjacent to any vertex in $V(P)\setminus\{v_2\}$. Then, the sequence $(a_2,c,a_1,v_1,z_1,\ldots,z_r)$ is a Z-sequence, in which each of $c,a_1$ and $v_2$ footprints exactly one vertex (notably, $v_2=z_r$ footprints only $u_2$), implying $\gr^Z(G)> \frac{n}{2}$, a contradiction. Now, we assume that $u_i$ is adjacent to a vertex in $V(P)\setminus\{v_i\}$, for each $i\in [2]$. If $u_1$ is adjacent to $z_1$, then the sequence $(a_2,c,a_1,v_1,z_1)$ is a Z-sequence in which each of $c,a_1$ and $z_1$ footprints exactly one vertex, which again yields $\gr^Z(G)> \frac{n}{2}$. By a similar argument, if $u_2$ is adjacent to $z_{r-1}$ we infer that  $\gr^Z(G)> \frac{n}{2}$. We may thus assume that $u_1$ is not adjacent to $z_1$, but it is adjacent to $z_2$, and that $u_2$ is not adjacent to $z_{r-1}$ but it is adjacent to $z_{r-2}$. Now, if $r>3$, consider the sequence $S=(a_2,c,a_1,v_1,z_1,z_2)$. It is easy to see that $S$ is a Z-sequence in which $c,a_1$ and $z_2$ footprint only one vertex. In an analogous  way as earlier we derive $\gr^Z(G)>\frac{n}{2}$, a contradiction.  If $r=3$, and $u_1u_2\notin E(G)$, then the sequence $(z_1,u_1,z_2,v_2,a_2,c)$ is a Z-sequence, and each of $z_2,v_2$ and $c$ footprints exactly one vertex, which again yields $\gr^Z(G)>\frac{n}{2}$. Finally, if $u_1u_2\in E(G)$, then the resulting graph $G$ is isomorphic to the necklace $N_{XY}$. 
 
The remaining possibility is that the only path between $v_1$ and $v_2$ is through $a_1$. Thus, $a_1v_1$ and $a_2v_2$ are cut-edges. Now, either $v_1$ already lies on a cycle, or there is a path between $v_1$ to a vertex $u$ in $G-a_1v_1$ that lies on a cycle, and let $u'$ be the vertex preceding $u$ on this path. By using Lemma~\ref{lem:Zlata}(i), we can extend the sequence $(a_1,c,a_2)$ (in which $c$ and $a_2$ footprint only one vertex) to a Z-dominating sequence in which all vertices, except $a_1$, footprint at most two vertices, and there is an additional vertex in the component $U$ of $G-u'u$ that contains $u$, which footprint only one vertex, unless $U$ is isomorphic to $X$. The possibility that $U$ is not isomorphic to $X$ leads to a contradiction. By Lemma~\ref{lem:Zlata}(ii), if $U$ is isomorphic to $X$, then there exists a dominating Z-sequence $S$, which starts with a vertex in $U$, and two vertices in $U$ of that sequence footprint only one vertex, while all other vertices, except the first vertex of $S$, footprint at most two vertices. This sequence can be chosen in such a way that contains $(v_1,a_1,c,a_2)$ as a consecutive subsequence, and so $c$ and $a_2$ are two additional vertices that footprint only one vertex (yielding in total at least four such vertices of $S$). This again implies the contradiction with the assumption $\gr^Z(G)> \frac{n}{2}$. 

\medskip

{\bf Case 3.} $G$ has a triangle, but is diamond-free.

Let $T$ be a triangle in $G$ with vertices $a_1,a_2,a_3$. Let $v_i$, where $i\in [3]$, be the neighbor of $a_i$ that does not belong to $T$. We start a sequence with $S=(a_1,a_2,a_3)$, which is clearly a Z-sequence, since $a_i$ footprints $v_i$. In addition, $a_2$ and $a_3$ footprint only one vertex. If $S$ is a dominating Z-sequence, then $\gr^Z(G)=3$ and $G=K_3\cp K_2$. Assume that $S$ is not a dominating Z-sequence, hence at least one of the vertices $v_i$ has at least one neighbor that is not in $T'=\{a_1,a_2,a_3,v_1,v_2,v_3\}$.  Assume that there is an edge between two vertices in $V(T')\setminus V(T)$, say $v_1v_2\in E(G)$. Then at least one of the vertices $v_1$ or $v_2$ has a neighbor outside $T'$. Thus, extending $S$ with that vertex (that is, with $v_1$ or $v_2$) yields a Z-sequence in which three vertices footprint exactly one vertex. In a similar way as in the previous cases, we derive that $\gr^Z(G)>\frac{n}{2}$. 

Suppose that there is a path between $v_1$ and $v_2$ that does not go through $a_1$ and $a_2$. Let $P:v_1=z_0,z_1,\ldots,z_r=v_2$ be a shortest such path. Let $x_i$ be the vertex in $N(v_i)\setminus V(P)$, which is not in $T$, for each $i\in [2]$. If $x_2$ is not adjacent to a vertex in $P$, then $(a_1,a_2,a_3,v_1,z_1,\ldots,z_r)$ is a Z-sequence in which each of $a_2,a_3$ and $v_2$ footprints exactly one vertex, yielding $\gr^Z(G)>\frac{n}{2}$. By an analogous argument we may also assume that $x_1$ is adjacent to a vertex in $P$. Suppose that $x_1$ is adjacent to $z_1$. Then $(a_1,a_2,a_3,v_1,z_1)$ is a Z-sequence in which $a_2,a_3$ and $z_1$ footprint exactly one vertex, and we are done. Hence, let $x_1$ be adjacent to $z_2$ ($x_1$ cannot be adjacent to $z_p$ for $p>2$, because this would yield a path from $v_1$ to $v_2$ shorter than $P$). Suppose first that $x_1\ne x_2$ (that is, the length $r$ of $P$ is greater than $2$). If $r>3$, then $(a_1,a_2,a_3,v_1,z_1,z_2)$ is a Z-sequence in which $a_2,a_3$ and $z_2$ each footprint only one vertex, and Lemma~\ref{lem:Zsequence} can be applied to infer $\gr^Z(G)>\frac{n}{2}$, a contradiction. Next, if $r=3$, then by symmetry we can also assume that $x_2$ is adjacent to $z_1$ (and $v_2$). If $x_1$ and $x_2$ are not adjacent, then the sequence $(a_1,a_2,a_3,v_1,x_1,v_2)$ is a Z-sequence, in which each of $a_2,a_3$, and $v_2$ footprints exactly one vertex (note that $v_2$ footprints $x_2$). On the other hand, if $x_1x_2\in E(G)$, then the sequence $S=(z_1,z_2,v_2,a_2,a_3)$ is a Z-sequence in which $z_1$ footprints $4$ vertices, $z_2$ footprints only $v_2$, $v_2$ footprints only $a_2$, and $a_3$ footprints only $v_3$, and we are done.  
Finally, if $r=2$, then $z_1$ is a common neighbor of $v_1$ and $v_2$. In this case, $(a_1,a_2,a_3,v_1,z_1)$ is a Z-sequence, and this time $z_1$ footprints only one vertex, unless $z_1$ is adjacent to $v_3$ (note that $z_1$ cannot be adjacent to $x_1=x_2$, since this would yield a diamond in $G$). The former possibility gives the desired bound, hence assume that $z_1v_3\in E(G)$. In addition, we may assume that $x_1$ is adjacent to $v_3$ for otherwise the roles of $z_1$ and $x_1$ can be reversed and we again get $\gr^Z(G)>\frac{n}{2}$. Now, this implies that $G$ has $8$ vertices, and is isomorphic to the graph $TK$ depicted on the left in Figure~\ref{fig:hamming}. 

The remaining case when there is no path between $v_1$ and $v_2$ that does not go through $a_1$ and $a_2$. By symmetry, we may assume that   $v_1a_1$ is a cut-edge. This situation can be dealt with in a similar way as in the last paragraph of Case 2. By using Lemma~\ref{lem:Zlata}(i), we infer that unless there is a subgraph $X$ that lies in a component $G-v_1a_1$ that contains $v_1$, we get $\grz(G)>\frac{n}{2}$. But, having $X$ as a subgraph, we can apply Lemma~\ref{lem:Zlata}(ii), and find a dominating Z-sequence $S$, which starts in $X$, two vertices of $X$ footprint only one vertex, and then at some point $S$ contains $(v_1,a_1,a_2,a_3)$ as a consecutive subsequence, yielding $\grz(G)>\frac{n}{2}$ again.

\medskip

{\bf Case 4.} $G$ is triangle-free. 

Let $C:v_1,\ldots, v_p,v_1$ be a shortest cycle in $G$. Since $G$ is triangle-free, $p\ge 4$, and obviously $C$ is an induced cycle. Moreover, if $p\ge 5$, then we claim that no two vertices of $C$ have a common neighbor outside $C$. Indeed, for $p=5$, two vertices with a common neighbor would imply that there exists either a triangle or a square in $G$, which is not possible since $C$ is a shortest cycle. When $p>5$, two vertices with a common neighbor would also imply that there is a cycle in $G$ with less vertices than $C$, a contradiction. Let us denote by $a_i$ the neighbor of $v_i$, which does not lie on $C$.

\medskip

\textbf{Case 4.A.} $\mathbf{p\ge 5}$.

If there is a vertex $a_i$ such that every path between $a_i$ and any vertex $a_j$, where $j\in [p]\setminus \{i\}$, passes through $C$, then $v_ia_i$ is a cut-edge. We derive, by using both statements of Lemma~\ref{lem:Zlata}, that there is a dominating Z-sequence, which starts in a subgraph isomorphic to $X$, and then eventually passes $(a_i,v_i,v_{i+1},\ldots, v_p,v_1,\ldots,v_{i-1})$ as a consecutive subsequence, in which there are at least four vertices that footprint only one vertex. This contradiction means that we may assume that for every $v_i$ there exists a vertex $v_j$ and a path between $v_i$ and $v_j$ that does not pass any vertex of $C$. Let $P: v_i=z_0,z_1=a_i,\ldots, z_{r-1}=a_j,z_r=v_j$ be a shortest such path over all $i,j \in [p]$. We observe that $z_k \neq a_k$, for any $k\in\{2,\ldots r-2\}$. 

Now, we again consider different subcases.

\textbf{Case 4.A.1.} $|V(P)| = 4$.

Note that this is possible only in two different cases. The first one is that $p=6$ and $v_i$ and $v_j$ are diametrical vertices of $C$. Therefore $(v_1,v_2,\ldots, v_6, a_i)$ is a Z-sequence in which $v_5,v_6$ and $a_i$ footprint exactly one vertex, a contradiction due to Lemma~\ref{lem:ZsequenceVec}. The second possibility is that $p=5$. Without loss of generality, let $v_i=v_1$ and $v_j=v_3$. If $a_1a_4 \notin E(G)$, then $(v_1,v_2,\ldots, v_5,a_1)$ is a Z-sequence in which $v_{4},v_5$ and $a_1$ footprint exactly one vertex. Suppose now that $a_1a_4 \in E(G)$. Then, if $a_3$ (resp.~$a_4$) is not adjacent to $a_5$ (resp.~$a_2$), the sequence $(v_1,v_2,\ldots, v_5,a_3 \textrm { (resp.} a_4))$ is Z-sequence and $v_{4},v_5$ and $a_3$ (resp.~$a_4$) footprint exactly one vertex. Now, if still $a_2$ and $a_5$ are not adjacent, the sequence $(v_1,v_2,\ldots, v_5,a_2)$ is a Z-sequence and $v_{4},v_5$ and $a_2$ footprint exactly one vertex. But if they are adjacent, we obtain the Petersen graph. 

\textbf{Case 4.A.2.} $|V(P)| = 5$.

Let $P:v_i,a_i, z_2, a_j, v_j$. Therefore $(v_1,v_2,\ldots, v_p, a_i,a_j)$ is a Z-sequence in which  $v_{p-1},v_p$ and $a_j$ footprint exactly one vertex. (Note, that we may choose $a_j$, since it cannot be adjacent to any of the vertices $a_k$, where $k \in [p]$, otherwise there would exist a shorter path between two vertices of a cycle $C$, contradicting the minimality of $P$.)

\textbf{Case 4.A.3.} $|V(P)| \ge 6$.
 
The sequence $(v_1,v_2,\ldots, v_p,z_1, z_2, \ldots, z_{r-2})$ is a Z-sequence in which $v_{p-1},v_p$ and $z_{r-2}$ footprint exactly one vertex. By using Lemma~\ref{lem:ZsequenceVec} we get $\gr^Z(G)>\frac{n}{2}$, a contradiction. 

\medskip 

\textbf{Case 4.B.} $\mathbf{p=4}$.

Suppose that there is a vertex $a_i$ (that is, the neighbor of $v_i$ outside $C$) such that $v_ia_i$ is a cut-edge in $G$. Similarly as in the previous case, we can apply both statements of Lemma~\ref{lem:Zlata}, and deduce that there exists a dominating Z-sequence $S$, which starts in a subgraph isomorphic to $X$, and then later passes $(a_i,v_i,v_{i+1},\ldots, v_p,v_1,\ldots,v_{i-1})$ as a consecutive subsequence; in $S$ there are at least four vertices that footprint only one vertex. This contradiction means that we may assume that $v_ia_i$ is not a cut-edge for all $i\in [4]$, and let $P$ be a shortest path between $v_i$ and $v_j$, among all pairs $\{i,j\}\subset [4]$, such that all internal vertices of $P$ are outside $C$.
Let $P:v_i,z_1=a_i,\ldots, z_{r-1}=a_j,v_j$, and observe that $z_k \neq a_k$, for any $k\in \{2,\ldots, r-2\}$. 

\textbf{Case 4.B.1.} $|V(P)| \ge 6$.

Clearly, in the sequence $S=(v_1,v_2, v_3, v_4,z_1, z_2, \ldots, z_{r-2})$ each of the vertices $v_{3}$, $v_4$ and $z_{r-2}$ footprints exactly one vertex, a contradiction. 

\textbf{Case 4.B.2.} $|V(P)| = 5$.

In that case, $S=(v_1,v_2,v_3,v_4,a_i,a_j)$ is a Z-sequence in which  $v_{p-1},v_p$ and $a_j$ footprint exactly one vertex. (Indeed, since $P$ is a shortest path with prescribed properties, there are no edges between $a_k$ and $a_\ell$ for distinct indices $k$ and $\ell$, which approves the choice of last two vertices in the sequence $S$, which is thus a Z-sequence.)

\textbf{Case 4.B.3.} $|V(P)| = 4$.

Let $v_1,v_2,v_3,v_4$ be the vertices of a cycle $C$ and let $a_1, a_2, a_3,a_4$ be their corresponding neighbors.  
(Note that $a_i \neq a_j$, for all $i,j \in [4]$, otherwise there would exist shorter path than $P$, a contradiction.)

First, consider the case when both vertices $a_i$ and $a_j$ are adjacent to some $a_m$ and $a_n$, respectively, where $m,n \in [p]\setminus \{i,j\}$. (Note, that $m \neq n$, otherwise we would have a triangle.) If $a_ma_n \notin E(G)$, then $(v_1,v_2,v_3,v_4, a_m)$ is a Z-sequence in which  $v_{3},v_4$ and $a_m$ footprint exactly one vertex. Otherwise, if $a_ma_n \in E(G)$, then $G$ is isomorphic to either $Q_3$ or $TQ_3$.

The second case is that some of the vertices $a_i$ or $a_j$ is not adjacent to any of $a_m$ and $a_n$, where $m,n \in [p]\setminus \{i,j\}$. Without loss of generality, let $a_i$ be such vertex. Hence, $a_i$ has a neighbor, different from $v_k$ and $a_k$ for all $k \in [4]$. Thus the sequence $(v_1,v_2,v_3,v_4,a_i)$ is a Z-sequence, and   $v_{3},v_4$ and $a_i$ footprint exactly one vertex, a contradiction. 

\textbf{Case 4.B.4.} $|V(P)| = 3$.

In this case, we may assume without loss of generality that $v_1$ and $v_3$ are open twins. Let $a$ be the common neighbor of $v_1$ and $v_3$, which is not in $C$. 
 
First, assume that $v_2$ and $v_4$ are also open twins, and let $b$ be the common neighbor of $v_2$ and $v_4$ outside $C$. Clearly, $a\ne b$, since $G$ has no triangles. Also, $a$ and $b$ are not adjacent, since then $G$ has $6$ vertices, and is isomorphic to $K_{3,3}$, which is a contradiction ($\grz(K_{3,3})=2$).  
If $a$ and $b$ have common neighbor $x$, then $x$ is a cut-vertex, has a neighbor $y\notin\{a,b\}$ such that $xy$ is a cut-edge. In a similar way as in the previous cases, we apply Lemma~\ref{lem:Zlata} to derive that $G\in \cal M$ (more precisely, $G$ is either $X_2$ or $X_3$).  
Suppose now, that $a$ and $b$ do not have common neighbor, and denote $c$ and $d$ their neighbors, respectively. If $cd \in E(G)$, then $(v_1,v_2,a,b,c)$ is a Z-sequence and each of $a$,$b$ and $c$ footprints exactly one vertex. On the other hand, let $cd \notin E(G)$, and first suppose that there exist a path between $c$ and $d$, which does not pass through the vertices $a$ and $b$. Let $P: c=z_0,z_1,\ldots,z_r=d$ be a shortest such path, and let us denote the neighbors of $c$ and $d$, which do not belong to $P$ by $c'$ and $d'$, respectively.  The rest of the proof of this situation goes along the same lines as in the second paragraph of Case 2 to obtain a third vertex that footprints only one vertex (noting that in the first part of the sequence there are two such vertices).  The only case in that proof, which does not lead to a contradiction is when $r=3$, $N(z_1)=N(d')$ and $N(z_2)=N(c')$ (which also implies $c'd' \in E(G)$). The resulting graph $G$ is isomorphic to the necklace $N_{XX}$; see Figure \ref{fig:necklace}. Finally, if $cd \notin E(G)$, but the only path between $c$ and $d$ goes through $a$ and $b$, then $ac$ and $bd$ are cut-edges. By using Lemma~\ref{lem:Zlata}, we infer that there is a dominating Z-sequence $S$, which starts in a subgraph isomorphic to $X$, two vertices of $X$ footprint only one vertex, and then at some point $S$ contains $(c,a,v_1,v_2,b)$ as a consecutive subsequence, and $b$ footprints only $d$, yielding $\grz(G)>\frac{n}{2}$, which is a contradiction.

The remaining case is that  $v_1$ and $v_3$ are open twins, while $v_2$ and $v_4$ are not open twins. In this case, vertices $v_1,\ldots, v_4$ and $a$ induce a $K_{2,3}$. Thus, $v_2,v_4$ and $a$ are in a symmetric role, and each of them has exactly one neighbor outside the mentioned $K_{2,3}$: let $x$ be the neighbor of $v_2$, $y$ the neighbor of $v_4$ and $z$ the neighbor of $a$. Since $G\ne K_{3,3}$ not all of the vertices $x,y,$ and $z$ coincide. It is possible that two of them coincide, yet this case can be then reduced to the previous one, where $v_2$ and $v_4$ were also open twins (possibly by renaming $a$ to $v_2$ or $v_4$). Hence, we may assume that $x,y$ and $z$ are pairwise distinct. 

If one of the edges $v_2x$, $v_4y$ or $az$ is a cut-edge, then the proof can be completed by using Lemma~\ref{lem:Zlata} as in several previous cases. We find that there exists a subgraph isomorphic to $X$, start a sequence in that subgraph, and continue it to reach the mentioned $K_{2,3}$ in such a way that at least three vertices footprint exactly one vertex, which yields a contradiction by Lemma~\ref{lem:ZsequenceVec}. We may therefore assume that none of the mentioned three edges is a cut-edge, and so there is a path between each of the pairs in $\{x,y,z\}$, which does not pass any vertex $v_i$. Assume that $P: x=z_0,z_1,\ldots,z_r=y$ be a shortest such path, and let us denote the neighbors of $x$ and $y$, which do not belong to $P$ by $x'$ and $y'$, respectively.  The rest of the proof of this situation goes along the same lines as in the second paragraph of Case 2 to obtain a third vertex that footprints only one vertex (in the initial part of the sequence $v_1,v_2,a,v_4)$, each of the vertices $a$ and $v_4$ footprints only one vertex). The only case in that proof, which does not lead to a contradiction is when $r=3$, $N(z_1)=N(x')$ and $N(z_2)=N(y')$ (which also implies $x'y' \in E(G)$). However, this implies that $az$ is a cut-edge, which was already considered above. The proof is complete.
\qed

\bigskip 

The consequence of Theorem~\ref{thm:Z-cubic} for the zero forcing number is immediate:

\begin{cor}
A connected, cubic, graph $G$ of order $n$ has $Z(G)=\frac{n}{2}$ if and only if $G\in {\cal M}'$ or $G$ is one of the graphs $N_{XX}, N_{XY}, N_{YY}$, $K_3\Box K_2$, $TK$, $Q_3$, $TQ_3$, or the Petersen graph.
\end{cor}

The following corollary for the Grundy domination number of cubic graphs follows from the fact that $\gr(G)\ge \grz(G)$ for all graphs $G$ with no isolated vertices, the value $\gr(K_{3,3})=3$, and by verifying the values of the Grundy domination numbers of $15$ extremal graphs from Theorem~\ref{thm:cubic}.
\begin{cor}
A connected, cubic, graph $G$ of order $n$ has $\gr(G)=\frac{n}{2}$ if and only if $G$ is one of the graphs $K_{3,3},Y_2, Y_3, N_{YY}, K_3\Box K_2$, $Q_3$, $TQ_3$, or the Petersen graph.
\end{cor}

\section*{Acknowledgments}

The first author was supported by the Ministry of Science of Slovenia under the grants P1-0297, J1-9109, J1-1693 and J1-2452.


\end{document}